\documentclass[draft,12pt]{article}
\usepackage{amsfonts}
\usepackage{amssymb}
\usepackage{amsbsy}
\usepackage{amsthm}
\usepackage{amsmath}
\usepackage{amscd}

\newtheorem{thm}{Theorem}[section]
\newtheorem{cor}[thm]{Corollary}
\newtheorem{lem}[thm]{Lemma}
\newtheorem{prop}[thm]{Proposition}
\theoremstyle{definition}
\newtheorem{defn}[thm]{Definition}
\theoremstyle{remark}
\newtheorem{rem}[thm]{Remark}
\newtheorem{ex}[thm]{\bf Example}

\newcommand{\bt}{\begin{thm}}
\newcommand{\et}{\end{thm}}
\newcommand{\bc}{\begin{cor}}
\newcommand{\ec}{\end{cor}}
\newcommand{\bl}{\begin{lem}}
\newcommand{\el}{\end{lem}}
\newcommand{\bp}{\begin{prop}}
\newcommand{\ep}{\end{prop}}
\newcommand{\bd}{\begin{defn}}
\newcommand{\ed}{\end{defn}}
\newcommand{\br}{\begin{rem}}
\newcommand{\er}{\end{rem}}
\newcommand{\bpr}{\begin{proof}}
\newcommand{\epr}{\end{proof}}
\newcommand{\bex}{\begin{ex}}
\newcommand{\eex}{\end{ex}}
\newcommand{\bcd}{\begin{CD}}
\newcommand{\ecd}{\end{CD}}

\newcommand{\bi}{\begin{itemize}}
\newcommand{\ei}{\end{itemize}}
\newcommand{\be}{\begin{enumerate}}
\newcommand{\ee}{\end{enumerate}}
\newcommand{\ba}{\begin{array}}
\newcommand{\ea}{\end{array}}
\newcommand{\beq}{\begin{equation}}
\newcommand{\eeq}{\end{equation}}
\newcommand{\beqa}{\begin{eqnarray}}
\newcommand{\eeqa}{\end{eqnarray}}
\newcommand{\bca}{\begin{cases}}
\newcommand{\eca}{\end{cases}}

\newcommand{\ds}{\displaystyle}
\newcommand{\ts}{\textstyle}

\newcommand{\N}{{\mathbb N}}

\newcommand{\R}{{\mathbb R}}
\newcommand{\C}{{\mathbb C}}
\newcommand{\T}{{\mathbb T}}
\newcommand{\D}{{\mathbb D}}

\newcommand{\U}{{\mathbb U}}
\newcommand{\OO}{{\mathbb O}}

\newcommand{\W}{{\mathbb W}}
\newcommand{\F}{{\mathbb F}}

\newcommand{\cB}{{\mathcal  B}}
\newcommand{\cC}{{\mathcal  C}}

\newcommand{\cH}{{\mathcal  H}}

\newcommand{\cS}{{\mathcal  S}}
\newcommand{\cT}{{\mathcal  T}}

\newcommand{\frB}{{\mathfrak B}}
\newcommand{\frC}{{\mathfrak C}}

\newcommand{\frM}{{\mathfrak M}}

\newcommand{\frP}{{\mathfrak P}}
\newcommand{\frS}{{\mathfrak S}}

\newcommand{\bs}{\boldsymbol}

\newcommand{\bsalpha}{{\boldsymbol \alpha}}

\newcommand{\bsgamma}{{\boldsymbol \gamma}}

\newcommand{\spn}{\mathrm{span}}
\newcommand{\supp}{\mathrm{supp}}

\newcommand{\Lim}{\mathrm{Lim\,}}
\newcommand{\re}{\mathrm{Re}}
\newcommand{\im}{\mathrm{Im}}

\begin{document}

\title{\bf Wall rational functions and Khrushchev's formula for orthogonal rational
functions\thanks{This work was partially realized during two stays of the second author at
the Norwegian University of Science and Technology (NTNU) financed respectively by Secretar\'{\i}a
de Estado de Universidades e Investigaci\'{o}n from the Ministry of Education and Science of Spain
and by the Department of Mathematical Sciences of NTNU.
The work of the second author was also partly supported by a research grant from the Ministry
of Education and Science of Spain, project code MTM2005-08648-C02-01, and by Project E-64
of Diputaci\'on General de Arag\'on (Spain).} }
\author{O. Nj\aa stad$^a$, L. Vel\'azquez$^b$ \\
\small{$^a$Norwegian University of Science and Technology, Norway} \\
\small{\texttt{njastad@math.ntnu.no}} \\
\small{$^b$University of Zaragoza, Spain} \\
\small{\texttt{velazque@unizar.es}}
}
\date{}
\maketitle

\vspace*{-5pt}

\begin{abstract}

We prove that the Nevalinna-Pick algorithm provides different homeomorphisms between
certain topological spaces of measures, analytic functions and sequences of complex
numbers. This algorithm also yields a continued fraction expansion of every Schur
function, whose approximants are identified. The approximants are quotients of rational
functions which can be understood as the rational analogs of the Wall polynomials. The
properties of these Wall rational functions and the corresponding approximants are
studied. The above results permit us to obtain a Khrushchev's formula for orthogonal
rational functions. An introduction to the convergence of the Wall approximants in the
indeterminate case is also presented.

\end{abstract}


\medskip

\noindent{\it (2000) AMS Mathematics Subject Classification}: 42C05.

\section{Introduction} \label{INT}

It is known that the Cayley transform provides a correspondence between Schur and
Carath\'{e}odory functions. Besides, the integral representation of Carath\'{e}odory functions
establishes a connection with finite positive Borel measures on the unit circle. On the
other hand, the Schur algorithm associates with any Schur function the so called Schur
parameters: a sequence in the open unit disk, the last point lying on the unit circle in
the case of a terminating sequence. Indeed, the set of these complex sequences, the set
of probability measures, the set of normalized Carath\'{e}odory functions and the set of
Schur functions become homeomorphic under suitable topologies.

The homeomorphism with the sequences of Schur parameters yields a bicontinuous
parametrization of the Schur functions or, alternatively, of the probability measures on
the unit circle. The study of such a parametrization is important, not only for the
theory of analytic functions, but also for the theory of continued fractions because
the Schur algorithm is equivalent to a continued fraction expansion of Schur functions
(hence, to a continued fraction for Carath\'{e}odory functions too). On the other hand, the
parametrization of measures on the unit circle becomes specially significant for the
associated orthogonal polynomials, since the Schur parameters are the coefficients of the
corresponding recurrence relation. The orthogonal polynomials on the unit circle also
provide the numerators and denominators of the approximants for the continued fraction
expansion of the related Carath\'{e}odory function. A similar role for the case of Schur
functions is played by the so called Wall polynomials, closely related to the orthogonal
polynomials too.

Therefore, the above homeomorphisms permit us to connect problems concerning measures,
orthogonal polynomials, continued fractions, analytic functions and complex sequences, so
that one can translate results or choose the best context to work. A remarkable example
of this is Krushchev's theory (see \cite{Kh01, Kh02}), which takes advantage of these
connections to reach deep and impressive results on the referred matters. A key result in
Khrushchev's theory is the so called Khrushchev's formula, obtained in \cite{Kh01}
starting from the analysis of the Wall polynomials. This formula can be understood as the
identification of the Schur functions of certain varying measures obtained by an
orthogonal polynomial modification of the orthogonality measure.

The Schur algorithm is a characterization of Schur functions based on an iteration which
evaluates each iterate at the origin. The Nevalinna-Pick algorithm, related to the
interpolation of Schur functions, generalizes this procedure evaluating each iterate at a
different point of the open unit disk. Like the Schur algorithm, the Nevalinna-Pick
generalization associates with any Schur function a similar sequence of parameters, but
depending now on the choice of the evaluation points. The Nevalinna-Pick algorithm is
also related to a rational generalization of the orthogonal polynomials on the unit
circle: the orthogonal rational functions with prescribed poles outside the unit circle.
It is known that these orthogonal rational functions are involved in alternative
continued fraction expansions of Carath\'{e}odory functions. However, the corresponding
continued fractions associated with Schur functions are not discussed in the literature.
The related approximants should have as numerators and denominators certain rational
functions depending on the evaluation points, which we will call Wall rational functions.

Finally we must comment a remakable new phenomenon of the Nevalinna-Pick algorithm which
does not appear in the Schur one: when the evaluation points approach to the unit circle
quickly enough, an indeterminate case can appear, i.e., different Schur functions can
have the same Nevalinna-Pick parameters. This causes important difficulties in the study
of the convergence of the corresponding continued fraction, which now can have different
limit points.

Once we have situated the context, we can understand the interest of our work, whose aims
are:

\vspace*{-10pt}

\bi

\item The analysis of the homeomorphisms related to the Nevalinna-Pick algorithm
(Section \ref{NP}).

\vspace*{-5pt}

\item The study of the Wall rational functions and the corresponding continued fraction
approximants of Schur functions (Section \ref{WR}).

\vspace*{-5pt}

\item The search for a Krushchev's formula for orthogonal rational functions (Section
\ref{KF}).

\vspace*{-5pt}

\item An introduction to the analysis of the limit points of continued fractions for Schur
functions in the indeterminate case (Section \ref{IC}).

\ei

\vspace*{-10pt}

We will follow Khrushchev's approach to the polynomial case given in \cite{Kh01,Kh02}. As
we will see, the approximants of a Schur function related to the Wall rational functions
are the Schur functions corresponding to the approximants of the continued fraction for
the related Carath\'{e}odory function. Hence, bearing in mind the homeomorphism between Schur
and Carath\'{e}odory functions, the convergence of the Schur continued fraction is equivalent
to the convergence of the Carath\'{e}odory continued fraction. Indeed, we will show that the
convergence of both continued fractions can be understood as a consequence of the
asymptotics of the Nevalinna-Pick parameters corresponding to the related approximants.
These convergence results are limited by the validity of the homeomorphism for the
Nevalinna-Pick parametrization of the Schur functions, which is ensured in the
determinate case.

The results about the Wall rational functions and the homeomorphisms related to the
Nevalinna-Pick algorithm will be the main tools to prove a Khrushchev's formula for
orthogonal rational functions. This will be the starting point of a ``rational
Khrushchev's theory" whose development will be given elsewhere. Nevertheless, a first
application of Khrushchev's formula will appear in the study of the indeterminate case.
The reason is that, contrary to the standard polynomial techniques, which usually can be
extended only to the determinate rational case, the rational generalization of
Khrushchev's formula always holds, providing an important tool for the study of the
indeterminate case. Nevertheless, our approach to the indeterminate case will be only
introductory, trying simply to show the variety of situations that can appear in the
convergence of the related continued fractions. A more complete study of the
indeterminate case deserves further investigations.

\section{Nevalinna-Pick homeomorphisms} \label{NP}

The results that we will prove here hold, not only for Schur functions on the unit disk
$\D=\{z\in\C:|z|<1\}$, but also for Schur functions on the upper half plane
$\U=\{z\in\C:\re\,z>0\}$, as can be seen using the Cayley transform. We will use a
unified notation to present simultaneously the results in both situations, and when we
want to distinguish between them we will write a left brace with the $\D$ case in the
first line and the $\U$ case in the second one. For instance, in what follows we will use
the notation
$$
\OO = \bca \D, \cr \U, \eca
\qquad
\partial\OO = \bca \T, \cr \overline\R, \eca
\qquad
\OO^e = \overline\C\setminus\overline\OO,
$$
where $\overline S$ is the closure in $\overline\C=\C\cup\{\infty\}$ of a subset
$S\subset\overline\C$.

Consider the transformations $\zeta_\alpha$, $\alpha\in\OO$, given by
$$
\zeta_\alpha=z_\alpha\frac{\varpi^*_\alpha}{\varpi_\alpha},
\qquad
z_\alpha=\bca -\frac{|\alpha|}{\alpha}, \smallskip \cr \frac{|1+\alpha^2|}{1+\alpha^2}, \eca
\quad
\varpi_\alpha(z)=\bca 1-\overline\alpha z, \cr z-\overline\alpha, \eca
\quad
\varpi^*_\alpha(z)=z-\alpha,
$$
where we understand that $z_\alpha=1$ for the particular value $\alpha=\alpha_0$ with
$$
\alpha_0 = \bca 0, \cr i. \eca
$$
$\zeta_\alpha$ is a homeomorphism of $\overline\C$ which maps $\OO$, $\partial\OO$ and
$\OO^e$ onto $\D$, $\T$ and $\overline\C\setminus\overline\D$ respectively.

A useful identity for $\zeta_\alpha$ is
\beq \label{ZETA-DIF}
\zeta_\alpha(t)-\zeta_\alpha(z) =
z_\alpha \frac{\varpi_\alpha(\alpha)\,\varpi^*_z(t)}{\varpi_\alpha(t)\,\varpi_\alpha(z)}.
\eeq
Besides, if we define the substar operation on complex functions by
$$
f_*(z)=\overline{f(\hat z)}, \qquad \hat z = \bca 1/\overline z, \cr \overline z, \eca
$$
then
\beq \label{ZETA*}
\qquad \zeta_{\alpha*}=\zeta_{\hat\alpha}=1/\zeta_\alpha.
\eeq

The sets that will be involved in the homeomorphisms are
$$
\kern-5pt
\ba{l} \frP = \hbox{ set of finite Borel measures on } \partial\OO,
\kern8pt \frP_0 = \{d\mu\in\frP : \mu_0=\int d\mu = 1\},
\medskip \cr
\frC = \{F\in\cH(\OO) : \re\,F(z)>0 \kern5pt \forall z\in\OO\},
\kern14pt \frC_\alpha = \{F\in\frC : F(\alpha)=1\},
\medskip \cr
\frB = \{f\in\cH(\OO) : |f(z)|\leq1 \kern5pt \forall z\in\OO\},
\medskip \cr
\frS = \left\{\bsgamma=(\gamma_n)_{n=0}^N : N\in\{0,1,\dots,\infty\},
\kern3pt \gamma_n \in \bca \D & $if $ n<N \cr \T & $if $ n=N<\infty \eca \right\},
\ea
$$
where $\alpha\in\OO$ and $\cH(S)$ is the set of analytic functions on the subset
$S\subset\C$. We will consider the topologies
$$
\begin{tabular}{l l l l}
set & \kern5pt topology & \kern5pt notation
\medskip \\
$\frP$ & \kern5pt $\ast$-weak convergence & \kern5pt $d\mu^k \stackrel{\ast}{\to}d\mu$
\medskip \\
$\frC,\frB$ & \kern5pt uniform convergence in compact subsets of $\OO$ & \kern5pt $f^k \rightrightarrows f$
\medskip \\
$\frS$ & \kern5pt pointwise convergence & \kern5pt $\bsgamma^k\to\bsgamma$
\end{tabular}
$$
The elements of $\frB$ and $\frC$ are called Schur and Carath\'{e}odory functions
respectively or, in short, S-functions and C-functions. We will assume that any Schur or
Carath\'{e}odory function $f$ is extended to $\partial\OO$ by
$$
f(z) =
\bca \lim_{r\uparrow1} f(rz), \cr \lim_{\epsilon\downarrow0} \re\,F(z+i\epsilon), \eca
\qquad \hbox{ a.e. } z\in\partial\OO,
$$
since it is known that such limits exits a.e. on $\partial\OO$.

The set of limit points of a sequence $(x^k)$ in a topological space will be denoted
$\Lim x^k$. We will use for the pointwise convergence in the space of complex sequences
the same notation as in the case of $\frS$. Concerning the convergence of an arbitrary
sequence $(f^k)$ of complex functions, the notation
$$
f^k \rightrightarrows f \text{ in } S
$$
means that $f^k$ converges uniformly to $f$ in compact subsets of $S\subset\C$.

For convenience, when using $\alpha_0$ as a subindex we will usually identify it with 0,
thus $\frC_0=\frC_{\alpha_0}$, $\zeta_0=\zeta_{\alpha_0}$, etc. In particular,
$$
\zeta_0(z)=\bca z, \smallskip \cr \ds \frac{z-i}{z+i}, \eca
$$
that is, $\zeta_0$ is the identity in $\D$ or the Cayley transform in $\U$.

\subsection{Measures, C-functions and S-functions}

Concerning the relation between measures and C-functions,
it is known that $\frC$ is homeomorphic to $\frP\times\R$ through
\beq \label{P-C}
\ba{c}
\mathop{\frP\times\R \xrightarrow{\kern23pt} \frC}
\limits_{\ds \hspace*{40pt} (d\mu,c) \to F(z;d\mu)+ic}
\smallskip \cr
F(z;d\mu) = \int D(t,z)\,d\mu(t), \qquad
\ds D(t,z) = \frac{\zeta_0(t)+\zeta_0(z)}{\zeta_0(t)-\zeta_0(z)}
= \bca
\frac{t+z}{t-z},
\smallskip \cr
\frac{1}{i} \frac{1+tz}{t-z}.
\eca
\ea
\eeq
In other words, $\frP$ is homeomorphic to the set of C-functions with the form $F(z;d\mu)$,
which are exactly the C-functions real valued at the origin since
$F(0;d\mu)=\mu_0$.
We have also the induced homeomorphism
\beq \label{P0-C0}
\cC_0 \colon \kern-15pt
\mathop{\frP_0 \xrightarrow{\kern20pt} \frC_0}
\limits_{\ds \hspace*{15pt} d\mu \to F(z;d\mu)}
\eeq
If $\mu'$ is the derivative of $d\mu$ with respect to the Lebesgue measure, it is known that
\beq \label{REF}
\re\,F(z;d\mu)=\mu'(z) \qquad \hbox{ a.e. } z\in\partial\OO.
\eeq

Some identities for $D(t,z)$ will be useful later. Let us start defining
$$
D_R(t,z)=\frac{1}{2}(D(t,z)+D_*(t,z)), \qquad D_I(t,z)=\frac{1}{2i}(D(t,z)-D_*(t,z)),
$$
where the substar operation on $D(\cdot,\cdot)$ is taken always on the first argument.
Then, $D_R(t,z)=\re\,D(t,z)$ and $D_I(t,z)=\im\,D(t,z)$ for $t\in\partial\OO$. Using
properties (\ref{ZETA-DIF}) and (\ref{ZETA*}) for $\zeta_0$ we find that
$$
D_*(t,z)=-D(t,\hat z)
$$
and
\beq \label{DD}
D(t,z)-D(t,\alpha) =
2 \frac{\varpi_0(t)\,\varpi^*_0(t)}{\varpi_0(\alpha_0)\,\varpi^*_z(t)}
\frac{\varpi^*_\alpha(z)}{\varpi^*_\alpha(t)}.
\eeq
Taking the substar operation with respect to $t$ on (\ref{DD}) and changing $z$ by
$\hat z$ we get
\beq \label{DD*}
D(t,z)+D_*(t,\alpha) =
2 \frac{\varpi_0(t)\,\varpi^*_0(t)}{\varpi_0(\alpha_0)\,\varpi^*_z(t)}
\frac{\varpi_\alpha(z)}{\varpi_\alpha(t)},
\eeq
which gives
\beq \label{DR}
D_R(t,z) = \frac{\varpi_z(z)}{\varpi_0(\alpha_0)}
\frac{\varpi_0(t)\,\varpi^*_0(t)}{\varpi_z(t)\,\varpi^*_z(t)}.
\eeq
In particular,
\beq \label{RED}
D_R(t,z) = \re\,D(t,z) = \frac{\varpi_z(z)}{\varpi_0(\alpha_0)}
\left|\frac{\varpi_0(t)}{\varpi_z(t)}\right|^2
=\bca
\frac{1-|z|^2}{|t-z|^2},
\smallskip \cr
\frac{\im\,z\,(1+t^2)}{|t-z|^2},
\eca
\qquad t\in\partial\OO.
\eeq

While the homeomorphism (\ref{P0-C0}) is the relevant one for the polynomial setting,
its generalization to $\frC_\alpha$ for any $\alpha\in\OO$ will be important for the
rational case. To understand this generalization notice that
$$
\mathop{\frP \xrightarrow{\kern20pt} \frP}
\limits_{\ds \hspace*{15pt} d\mu \to
\text{\footnotesize$\frac{d\mu(\cdot)}{D_R(\cdot,\alpha)}$}}
$$
is a homeomorphism since $D_R(t,\alpha)$ is positive and continuous for all $t\in\partial\OO$.
Composing it with (\ref{P-C}) shows that
$$
\mathop{\frP\times\R \xrightarrow{\kern35pt} \frC}
\limits_{\ds \hspace*{55pt} (d\mu,c) \to
F\kern-2pt\left(z;\text{\scriptsize$\frac{d\mu(\cdot)}{D_R(\cdot,\alpha)}$}\right)+ic}
$$
is a homeomorphism too. It is straightforward to see that this homemomorphism induces the
following one
\beq \label{P0-Ca}
\cC_\alpha \colon \kern-18pt
\mathop{\frP_0 \xrightarrow{\kern20pt} \frC_\alpha}
\limits_{\ds \hspace*{18pt} d\mu \to F_\alpha(z;d\mu)}
\eeq
where $F_\alpha(z;d\mu)$ is defined for any $d\mu\in\frP$ and any $\alpha\in\OO$ by
\vskip10pt
$
F_\alpha(z;d\mu) =
F\kern-2pt\left(z;\frac{d\mu(\cdot)}{D_R(\cdot,\alpha)}\right) + ic_\alpha(d\mu),
\qquad c_\alpha(d\mu)=-\int \frac{D_I(t,\alpha)}{D_R(t,\alpha)}\,d\mu(t).
$
\vskip10pt
\noindent We will say that $F_\alpha(z;d\mu)$ is the $\alpha$-C-function of $d\mu$.
$F_\alpha(\alpha;d\mu)=\mu_0$, thus, a C-function has the form
$F_\alpha(z;d\mu)$ for some $d\mu\in\frP$ iff it is real at $\alpha$, and the set
of these C-functions is homeomorphic to $\frP$.

A stronger convergence property than the one given by the homeomorphism (\ref{P0-Ca})
holds. To prove it we will use an explicit relation between $F_\alpha(z;d\mu)$ and
$F(z;d\mu)$. Although such a relation was obtained in \cite[Lemmas 6.2.2 and
6.2.3]{BGHN99}, we present here a more concise proof which, at the same time, unifies the
discussion for measures on the unit circle and the real line.

\bp \label{F-Fa}

For any $d\mu\in\frP$ and any $\alpha\in\OO$,
$$
F(z;d\mu) = D_R(z,\alpha)\,F_\alpha(z;d\mu) - i\mu_0D_I(z,\alpha).
$$

\ep

\bpr
From (\ref{DD}) and (\ref{DD*}) we find that
$$
\ba{l}
D(t,z)-iD_I(t,\alpha) =
\frac{1}{2} (D(t,z)-D(t,\alpha)) +
\frac{1}{2}(D(t,z)-D_*(t,\alpha))=
\medskip \cr \ds \kern98pt
= \frac{\varpi_0(t)\,\varpi^*_0(t)}{\varpi_0(\alpha_0)\,\varpi^*_z(t)}
\left( \frac{\varpi^*_\alpha(z)}{\varpi^*_\alpha(t)} +
\frac{\varpi_\alpha(z)}{\varpi_\alpha(t)} \right),
\ea
$$
which combined with (\ref{DR}) gives
$$
\frac{D(t,z)}{D_R(t,\alpha)}-i\frac{D_I(t,\alpha)}{D_R(t,\alpha)} =
\frac{\varpi_\alpha(t)\,\varpi^*_\alpha(z)+\varpi^*_\alpha(t)\,\varpi_\alpha(z)}
{\varpi_\alpha(\alpha)\,(t-z)}.
$$
The above function, as well as $D(t,z)$, are antisymmetric under the exchange of $t$ and
$z$. Hence,
$$
\frac{D(t,z)}{D_R(t,\alpha)} - i\frac{D_I(t,\alpha)}{D_R(t,\alpha)} =
\frac{D(t,z)}{D_R(z,\alpha)} + i\frac{D_I(z,\alpha)}{D_R(z,\alpha)},
$$
which, integrated with respect to $d\mu(t)$, finally yields the result.
\epr

The above relation permits us to obtain a convergence property for sequences of
C-functions normalized at different points. Notice that, as a consequence of the maximum
modulus principle, if $F^k,F\in\cH(\OO)$ and $F^k \rightrightarrows F$ in $\OO \setminus
K$, $K$ a compact subset of $\OO$, then $F^k \rightrightarrows F$ in $\OO$. If, besides,
$F^k\in\frC_{\alpha_k}$ with $\alpha_k\in\OO$ such that $\alpha_k\to\alpha\in\OO$, then
$F\in\frC_\alpha$, so we can suppose in this situation that
$F^k(z)=F_{\alpha_k}(z;d\mu^k)$ and $F(z)=F_\alpha(z;d\mu)$ for some probability measures
$d\mu^k$, $d\mu$.

\bt \label{CONV-Fa}

Let $(\alpha_k)$ be a sequence in $\OO$ and $(d\mu^k)$ a sequence in $\frP$. If
$\alpha_k\to\alpha\in\OO$, then
$$
F_{\alpha_k}(z;d\mu^k) \rightrightarrows F_\alpha(z;d\mu)
\kern5pt \Longleftrightarrow \kern5pt
d\mu^k\kern-3pt\stackrel{*}{\to}d\mu.
$$

\et

\bpr
It suffices to prove $ F_{\alpha_k}(z;d\mu^k) \rightrightarrows F_\alpha(z;d\mu)
\Leftrightarrow F(z;d\mu^k) \rightrightarrows F(z;d\mu). $ If $F_{\alpha_k}(z;d\mu^k)
\rightrightarrows F_\alpha(z;d\mu)$, then
$\mu^k_0=F_{\alpha_k}(\alpha_k;d\mu^k)\to\mu_0=F_\alpha(\alpha;d\mu)$. From this result,
Proposition \ref{F-Fa}, and the fact that $D_R(\cdot,\alpha_k) \rightrightarrows
D_R(\cdot,\alpha)$, $D_I(\cdot,\alpha_k) \rightrightarrows D_I(\cdot,\alpha)$ in
$\OO\setminus\{\alpha\}$, we conclude that $F(z;d\mu^k) \rightrightarrows F(z;d\mu)$ in
$\OO\setminus\{\alpha\}$, so it holds in $\OO$ too. A similar reasoning proves the
opposite implication, bearing in mind that $1/D_R(\cdot,\alpha_k) \rightrightarrows
1/D_R(\cdot,\alpha)$ in $\OO\setminus\{\alpha_0\}$.
\epr

As for the connection with Schur functions, the relations
$$
f = \frac{1}{\zeta_\alpha}\frac{F-1}{F+1},
\qquad
F = \frac{1+\zeta_\alpha f}{1-\zeta_\alpha f},
$$
define a one to one mapping between C-functions $F\in\frC_\alpha$ and S-functions
$f\in\frB$.
Moreover, for any $\alpha\in\OO$, the bijection
\beq \label{Ca-B}
\cB_\alpha \colon \kern0pt
\mathop{\frC_\alpha \xrightarrow{\kern0pt} \frB}
\limits_{\ds \hspace*{0pt} F \to f}
\eeq
is also a homeomorphism, as the following more general property shows.

\bt \label{CONV-f}

Let $(\alpha_k)$, $(\beta_k)$ be two sequences compactly included in $\OO$,
$F^k\in\frC_{\alpha_k}$, $G^k\in\frC_{\beta_k}$ and $f^k=\cB_{\alpha_k}(F^k)$,
$g^k=\cB_{\beta_k}(G^k)$. If $\alpha_k-\beta_k\to0$, then
$$
F^k - G^k \rightrightarrows 0
\kern5pt \Longleftrightarrow \kern5pt
f^k - g^k \rightrightarrows 0.
$$

\et

\bpr
The result follows easily from the identities
$$
\ba{l}
\ds f^k-g^k = \frac{1}{\zeta_{\alpha_k}}
\left[ (\zeta_{\beta_k}-\zeta_{\alpha_k})g^k + 2\frac{F^k-G^k}{(F^k+1)(G^k+1)} \right],
\medskip \cr
\ds F^k-G^k = 2\frac{(\zeta_{\alpha_k}-\zeta_{\beta_k})f^k+\zeta_{\beta_k}(f^k-g^k)}
{(1-\zeta_{\alpha_k}f^k)(1-\zeta_{\beta_k}g^k)},
\ea
$$
which give in $\OO$ the inequalities
$$
|f^k-g^k| \leq
\frac{|\zeta_{\alpha_k}-\zeta_{\beta_k}|+2|F^k-G^k|}{|\zeta_{\alpha_k}|},
\quad
|F^k-G^k| \leq 2
\frac{|\zeta_{\alpha_k}-\zeta_{\beta_k}|+|f^k-g^k|}{(1-|\zeta_{\alpha_k}|)(1-|\zeta_{\beta_k}|)}.
$$
When $(\alpha_k)$, $(\beta_k)$ are compactly included in $\OO$ and
$\alpha_k-\beta_k\to0$, the above inequalities prove that $f^k-g^k\rightrightarrows0$
implies $F^k-G^k\rightrightarrows0$, while $F^k-G^k\rightrightarrows0$ implies
$f^k-g^k\rightrightarrows0$ in $\OO\setminus\Lim\alpha_k$, thus in $\OO$ because
$\Lim\alpha_k$ is a compact subset of $\OO$.
\epr

Given $d\mu\in\frP_0$, the S-function $f_\alpha(z;d\mu)=\cB_\alpha(F_\alpha(z;d\mu))$,
will be called the $\alpha$-S-function of $d\mu$. The relation between $F_\alpha(z;d\mu)$
and $F(z;d\mu)$ provides an explicit expression of the $\alpha$-S-function
$f_\alpha(z;d\mu)$ of $d\mu$ in terms of its $\alpha_0$-S-function
$f(z;d\mu)=f_{\alpha_0}(z;d\mu)$.

\bp \label{f-fa}

Let $d\mu\in\frP_0$ and $\alpha\in\OO$. Denoting $f(z)=f(z;d\mu)$ and
$f_\alpha(z)=f_\alpha(z;d\mu)$,
$$
f_\alpha = -\frac{\zeta_0(\alpha)}{|\zeta_0(\alpha)|}
\frac{f-\overline{\zeta_0(\alpha)}}{1-\zeta_0(\alpha)f},
\qquad
f = -\frac{|\zeta_0(\alpha)|}{\zeta_0(\alpha)}
\frac{f_\alpha-|\zeta_0(\alpha)|}{1-|\zeta_0(\alpha)|f_\alpha}.
$$

\ep

\bpr
From Proposition \ref{F-Fa} we find that
$$
f_\alpha = \frac{1}{\zeta_\alpha}
\frac{(1-D_{\alpha*})+(1+D_{\alpha*})\,\zeta_0f}
{(1+D_\alpha)+(1-D_\alpha)\,\zeta_0f},
\qquad
D_\alpha(z)=D(z,\alpha).
$$
Besides, a direct calculation using the properties of $\zeta_0$ gives
$$
\frac{1-D_{\alpha*}}{1+D_\alpha} = |\zeta_0(\alpha)|\,\zeta_\alpha,
\qquad
\frac{1+D_{\alpha*}}{1-D_\alpha} = \frac{\zeta_\alpha}{|\zeta_0(\alpha)|},
\qquad
\frac{1+D_\alpha}{1-D_\alpha} = -\frac{\zeta_0}{\zeta_0(\alpha)}.
$$
From the first two identities we obtain
$$
f_\alpha =
\frac{|\zeta_0(\alpha)|(1+D_\alpha)+|\zeta_0(\alpha)|^{-1}(1-D_\alpha)\,\zeta_0f}
{(1+D_\alpha)+(1-D_\alpha)\,\zeta_0f},
$$
and, then, the last of the three identities yields the result.
\epr

In what follows, we will refer to $F(z;d\mu)$ and $f(z;d\mu)$ as the C-function and
S-function of $d\mu$ respectively.

Proposition \ref{f-fa} can be combined with Theorem \ref{CONV-f} to give the following
general equivalences.

\bt \label{CONV-gral}

Let $(\alpha_k)$, $(\beta_k)$ be two sequences compactly included in $\OO$ and
$(d\mu^k)$, $(d\nu^k)$ two sequences in $\frP_0$. If $\alpha_k-\beta_k\to0$, then
$$
\ba{l}
F_{\alpha_k}(z;d\mu^k)-F_{\beta_k}(z;d\nu^k)\rightrightarrows0
\kern5pt \Longleftrightarrow \kern5pt
f_{\alpha_k}(z;d\mu^k)-f_{\beta_k}(z;d\nu^k)\rightrightarrows0
\kern5pt \Longleftrightarrow \kern5pt
\smallskip \cr
\kern160pt \Longleftrightarrow \kern5pt
d\mu^k\kern-1pt-d\nu^k\stackrel{*}{\to}0.
\ea
$$

\et

\bpr
Suppose that $(\alpha_k)$, $(\beta_k)$ are compactly included in $\OO$ and
$\alpha_k-\beta_k\to0$. Then, Theorem \ref{CONV-f} ensures the first equivalence. With
the help of Proposition \ref{f-fa} we find that
$f_{\alpha_k}(z;d\mu^k)-f_{\beta_k}(z;d\nu^k)\rightrightarrows0$ iff
$f(z;d\mu^k)-f(z;d\nu^k)\rightrightarrows0$. Applying again Theorem \ref{CONV-f} we
conclude that $f(z;d\mu^k)-f(z;d\nu^k)\rightrightarrows0$ iff
$F(z;d\mu^k)-F(z;d\nu^k)\rightrightarrows0$. This last condition is equivalent to
$d\mu^k-d\nu^k\stackrel{*}{\to}0$ because
$F(z;d\mu^k)-F(z;d\nu^k)=F(z;d\mu^k\kern-2pt-d\nu^k)$.
\epr

\bex \label{LEBESGUE}
Let us define for any $\alpha\in\OO$ the measure
$$
dm_\alpha(t) =
\frac{\varpi_\alpha(\alpha)}{\varpi_\alpha(t)\,\varpi^*_\alpha(t)} \frac{dt}{2 \pi i} =
\begin{cases}
\frac{1-|\alpha|^2}{|t-\alpha|^2} \frac{dt}{2 \pi i t},
\smallskip \cr
\frac{\im\,\alpha}{|t-\alpha|^2} \frac{dt}{\pi},
\end{cases}
\qquad t\in\partial\OO.
$$

In particular, $dm=dm_{\alpha_0}$ is the Lebesgue measure
in $\T$ or its Cayley transform in $\R$, i.e.,
$$
dm(t) =
\begin{cases}
\frac{dt}{2 \pi i t},
\smallskip \cr
\frac{dt}{\pi(1+t^2)},
\end{cases}
\qquad t\in\partial\OO.
$$
Therefore, $F(z;dm)=1$, so $f(z;dm)=0$ and, from Proposition \ref{f-fa},
$$
f_\alpha(z;dm)=|\zeta_0(\alpha)|,
\qquad
F_\alpha(z;dm) =
\frac{1+|\zeta_0(\alpha)|\,\zeta_\alpha(z)}{1-|\zeta_0(\alpha)|\,\zeta_\alpha(z)}.
$$

In the general case, from (\ref{DR}) we get $dm_\alpha(t)=D_R(t,\alpha)\,dm(t)$,
thus we have the equality
\vskip10pt
\centerline{
$
F\left(z;\frac{dm_\alpha(\cdot)}{D_R(\cdot,\alpha)}\right) =
F(z;dm)=1,
$
}
\vskip10pt
\noindent which gives
\vskip10pt
$
\int dm_\alpha =
\re\,F\left(z;\frac{dm_\alpha(\cdot)}{D_R(\cdot,\alpha)}\right)
= 1,
\qquad
c_\alpha(dm_\alpha) =
-\im\,F\left(z;\frac{dm_\alpha(\cdot)}{D_R(\cdot,\alpha)}\right)
= 0.
$
\vskip10pt
\noindent Hence, $dm_\alpha\in\frP_0$ and
$$
F_\alpha(z;dm_\alpha)=1, \qquad f_\alpha(z;dm_\alpha)=0.
$$
Besides, Proposition \ref{f-fa} implies that
$$
f(z;dm_\alpha)=\overline{\zeta_0(\alpha)},
\qquad
F(z;dm_\alpha) =
\frac{1+\overline{\zeta_0(\alpha)}\,\zeta_\alpha(z)}{1-\overline{\zeta_0(\alpha)}\,\zeta_\alpha(z)}.
$$
This shows that the homeomorphism $\cB_0\cC_0$ between $\frP_0$ and $\frB$ establishes a
one to one correspondence between the set of measures $\{dm_\alpha : \alpha\in\OO\}$ and
the set of constant functions with values in $\D$.

The rest of constant S-functions are the constant unimodular ones, which the
homeomorphism $\cB_0\cC_0$ puts in one to one correspondence with the set
$\{\delta_\tau(t)=\delta(t-\tau)\,dt : \tau\in\partial\OO\}$ of Dirac measures, since
$$
F(z;\delta_\tau)=\frac{1+\overline{\zeta_0(\tau)}\,\zeta_0(z)}{1-\overline{\zeta_0(\tau)}\,\zeta_0(z)},
\qquad
f(z;\delta_\tau)=\overline{\zeta_0(\tau)}.
$$
The fact that
$f(z;dm_\alpha)=\overline{\zeta_0(\alpha)}
\underset{\alpha\to\tau}{\rightrightarrows}
f(z;\delta_\tau)=\overline{\zeta_0(\tau)}
$
implies $dm_\alpha \stackrel{*}{\underset{\alpha\to\tau}{\longrightarrow}} \delta_\tau$.
{\huge$\hfill\centerdot$}

\eex

The previous results deal only with the case of sequences $\bsalpha=(\alpha_n)$ compactly
included in $\OO$. If $\bsalpha$ is in $\OO$ but not compactly included there, a
subsequence $(\alpha_{n_j})_j$ must exist such that
$\Lim_{\!j}\,\alpha_{n_j}\subset\partial\OO$, i.e.,
$|\zeta_0(\alpha_{n_j})|\stackrel{j}{\to}1$. Concerning this situation we have the
following strong convergence result.

\bt \label{CONV-borde}

If $(\alpha_k)$ is a sequence in $\OO$ such that $\Lim\alpha_k\subset\partial\OO$, and
$(d\mu^k)$ is a sequence in $\frP_0$,
$$
1\notin\Lim f_{\alpha_k}(z;d\mu^k)
\kern5pt \Longrightarrow \kern5pt
d\mu^k-dm_{\alpha_k}\stackrel{*}{\to}0,
\kern7pt
\Lim d\mu^k=\{\delta_\tau:\tau\in\Lim\alpha_k\}.
$$

\et

\bpr
The relation between $g^k(z)=f_{\alpha_k}(z;d\mu^k)$ and $f^k(z)=f(z;d\mu^k)$ given by
Proposition \ref{f-fa} can be written as
$$
(1-g^k) \left(f^k-\frac{1}{\zeta_0(\alpha_k)}\right) +
\frac{1-|\zeta_0(\alpha_k)|}{|\zeta_0(\alpha_k)|}
\left(f^k+\frac{|\zeta_0(\alpha_k)|}{\zeta_0(\alpha_k)}\right)=0.
$$
The condition $\Lim\alpha_k\subset\partial\OO$ is equivalent to
$|\zeta_0(\alpha_k)|\to1$, so
$$
(1-g^k) \left(f^k-\overline{\zeta_0(\alpha_k)}\right) \rightrightarrows 0.
$$
The fact that $1\notin\Lim g^k$ forces
$f^k-\overline{\zeta_0(\alpha_k)}\rightrightarrows0$, which, in view of Example
\ref{LEBESGUE}, means that $f(z;d\mu^k)-f(z;dm_{\alpha_k})\rightrightarrows0$. Then, the
results follow from Theorem \ref{CONV-gral} and the last comment of Example
\ref{LEBESGUE}.
\epr

\subsection{The Nevalinna-Pick algorithm and the orthogonal rational functions}

The Nevalinna-Pick algorithm comes from the fact that the transformation
$$
f \to \frac{1}{\zeta_\alpha}\frac{f-f(\alpha)}{1-\overline{f(\alpha)}f}
$$
maps the interior $\frB^0=\{f\in\cH(\OO) : |f(z)|<1 \kern5pt \forall z\in\OO\}$ of $\frB$
on $\frB$ for any $\alpha\in\OO$. Given a sequence $\bsalpha=(\alpha_n)$ in $\OO$, this
algorithm associates with any $f\in\frB$ a finite or infinite sequence $(f_n)$ in $\frB$
defined by
\beq \label{NP-ITERATES}
\ba{l}
f_0=f,
\smallskip \cr \ds
f_{n+1} = \frac{1}{\zeta_{n+1}}\frac{f_n-\gamma_n}{1-\overline\gamma_nf_n},
\qquad \gamma_n=f_n(\alpha_{n+1}), \qquad \zeta_n=\zeta_{\alpha_n},
\qquad n\geq0,
\ea
\eeq
so that the sequence terminates at $f_N$ iff $f_N\in\frB\setminus\frB^0$,
which holds iff $f$ is, up to a unimodular factor, a finite Blaschke product
$\zeta_{\beta_1}\zeta_{\beta_2}\cdots\zeta_{\beta_N}$ with $\beta_k\in\OO$ for all $k$.
We will say that $(f_n)$ are the $\bsalpha$-iterates of $f$ and $\bsgamma=(\gamma_n)$ the
$\bsalpha$-parameters of $f$. Notice that $(f_n,f_{n+1},\dots)$ and
$(\gamma_n,\gamma_{n+1},\dots)$ are the iterates and parameters of $f_n$ associated
with the sequence $(\alpha_{n+1},\alpha_{n+2},\dots)$. From the relation between
$f_n$ and $f_{n+1}$ we easily obtain
\beq \label{ID}
(1-\overline\gamma_nf_n)(1+\overline\gamma_n\zeta_{n+1}f_{n+1})=1-|\gamma_n|^2,
\eeq
an identity which will be useful later.

The maximum modulus principle implies that
$\frB\setminus\frB^0$ is the set of constant unimodular functions. Therefore,
$\bsgamma\in\frS$. Indeed, the map
\beq \label{B-S}
\cT_\bsalpha \colon \kern0pt
\mathop{\frB \xrightarrow{\kern0pt} \frS}
\limits_{\ds \hspace*{0pt} f \to \bsgamma}
\eeq
is continuous for any sequence $\bsalpha$ in $\OO$, as follows from the following theorem,
which states a stronger result.

\bt \label{CONV-gamma}

Let $(\bsalpha^k)$ be a sequence of sequences in $\OO$, $\bsalpha$ a sequence in $\OO$,
$(f^k)$ a sequence in $\frB$, and $f\in\frB$. Then,
$$
\bsalpha^k \to \bsalpha,
\kern5pt f^k \rightrightarrows f
\kern5pt \Longrightarrow \kern5pt
\bsgamma^k \to \bsgamma, \kern5pt
f^k_n \stackrel{k}{\rightrightarrows} f_n \kern5pt \forall n,
$$
where $\bsgamma^k$ and $(f^k_n)_n$ are the $\bsalpha^k$-parameters and
$\bsalpha^k$-iterates of $f^k$, while $\bsgamma$ and $(f_n)$ are the
$\bsalpha$-parameters and $\bsalpha$-iterates of $f$, respectively.

\et

\bpr
Let $\bsgamma=(\gamma_n)_{n=0}^N$. We must prove that, for each $n \leq N$, $f^k_n$,
$\gamma^k_n$ exist for big enough $k$ and $f^k_n \stackrel{k}{\rightrightarrows} f_n$,
$\gamma^k_n \stackrel{k}{\to} \gamma_n$. Let us proceed by induction. First, $f^k_0=f^k$
exists for any $k$ and $f^k_0 \rightrightarrows f_0=f$ from the hypothesis. Now, given $n
\leq N$, suppose that $f^k_n$ exist for big enough $k$ and
$f^k_n\stackrel{k}{\rightrightarrows}f_n$. Then, $\gamma^k_n=f^k_n(\alpha^k_{n+1})$
exists for the same values of $k$ and
$\gamma^k_n\stackrel{k}{\to}\gamma_n=f_n(\alpha_{n+1})$. If $n=N$, there is nothing more
to prove. Otherwise, $f_n\in\frB^0$, thus $f^k_n\in\frB^0$ for big enough $k$ because
$f^k_n \stackrel{k}{\rightrightarrows} f_n$. In consequence, $f^k_{n+1}$ exists for such
values of $k$. Moreover, denoting $\zeta^k_n=\zeta_{\alpha^k_n}$,
$$
\ba{l}
\ds f^k_{n+1}-f_{n+1} =
\zeta^k_{n+1}f^k_{n+1} \left( \frac{1}{\zeta^k_{n+1}}-\frac{1}{\zeta_{n+1}} \right) +
\frac{1}{\zeta_{n+1}} \left( \zeta^k_{n+1}f^k_{n+1}-\zeta_{n+1}f_{n+1} \right) =
\smallskip \cr
\ds = \frac{1}{\zeta_{n+1}} \bigg[ (\zeta_{n+1}-\zeta^k_{n+1})f^k_{n+1}+
\smallskip \cr \kern60pt
\ds +\frac{(f^k_n-f_n)+(\gamma_n-\gamma^k_n)+(\overline\gamma^k_n-\overline\gamma_n)f^k_nf_n
+\gamma^k_n\overline\gamma_nf_n-\overline\gamma^k_n\gamma_nf^k_n}
{(1-\overline\gamma^k_nf^k_n)(1-\overline\gamma_nf_n)} \bigg],
\ea
$$
thus,
$$
|f^k_{n+1}-f_{n+1}| \leq \frac{1}{|\zeta_{n+1}|}
\left[ |\zeta^k_n-\zeta_n|
+\frac{2|f^k_n-f_n|+4|\gamma^k_n-\gamma_n|}{(1-|\gamma_n|)(1-|\gamma^k_n|)} \right],
$$
proving that $f^k_{n+1}\stackrel{k}{\rightrightarrows}f_{n+1}$ in
$\OO\setminus\{\alpha_{n+1}\}$, hence, in $\OO$.
\epr

Summarizing, given $\alpha\in\OO$ and an arbitrary sequence $\bsalpha=(\alpha_n)$ in $\OO$,
we have the following chain
$$
\bcd
\frP_0   @>\cC_\alpha>>   \frC_\alpha   @>\cB_\alpha>>   \frB   @>\cT_\bsalpha>>   \frS
\ecd
$$
the first two maps being homeomorphisms and the last one being continuous. For the choice
$\alpha=\alpha_0$, the above diagram can be closed to a commutative one. This result is a
consequence of the relation between $S$-functions and orthogonal rational functions.

Given a measure $d\mu\in\frP_0$ and a sequence $\bsalpha=(\alpha_n)$ in $\OO$, we can
consider the orthonormalization in $L^2(d\mu)$ of the Blaschke products $(B_n)$ given by
$$
\ba{l}
B_0=1
\smallskip \cr
B_n=\zeta_1\zeta_2\cdots\zeta_n, \kern10pt n\geq1.
\ea
$$
The result are the so called orthogonal rational functions $(\Phi_n)$ associated with
$d\mu$ and $\bsalpha$. Under a suitable normalization, they satisfy the recurrence relation
(see \cite[Theorem 4.1.3]{BGHN99})
\beq \label{RR0}
\ba{l}
\Phi_0=1,
\medskip \cr \ds
\begin{pmatrix} \Phi_n \cr \Phi^*_n \end{pmatrix} =
e_n {\varpi_{n-1}\over\varpi_n}
\begin{pmatrix} 1 & \overline\Lambda_n \cr \Lambda_n & 1 \end{pmatrix}
\begin{pmatrix} z_n \overline z_{n-1} \zeta_{n-1} \Phi_{n-1} \cr \Phi^*_{n-1} \end{pmatrix},
\qquad n\geq1,
\medskip \cr \ds
\Phi^*_n=B_n\Phi_{n*}, \qquad
\Lambda_n\in\D,\qquad
e_n=\ts\sqrt{{\varpi_n(\alpha_n)\over\varpi_{n-1}(\alpha_{n-1})}{1 \over 1-|\Lambda_n|^2}},
\ea
\eeq
where, for convenience, when $\alpha_n$ is a subindex it is denoted by $n$. In what follows,
when referring to orthogonal rational functions we will suppose that they are normalized so that
(\ref{RR0}) holds.
For our purposes, a more appropriate form of the above recurrence is in terms of the
functions $\hat\Phi_n=\overline z_n\Phi_n$ and the parameters
$\lambda_n=-z_{n+1}\Lambda_{n+1}$, i.e.,
\beq \label{RR}
\begin{pmatrix} \hat\Phi_n \cr \Phi_n^* \end{pmatrix} =
e_n {\varpi_{n-1}\over\varpi_n} \, T_{n-1}
\begin{pmatrix} \hat\Phi_{n-1} \cr \Phi_{n-1}^* \end{pmatrix},
\quad
T_n=\begin{pmatrix} \zeta_n & -\overline\lambda_n \cr -\lambda_n\zeta_n & 1 \end{pmatrix}.
\eeq
Notice that $\lambda_n\in\D$ is given by
$\lambda_n=-z_{n+1}\overline{\Phi_{n+1}(\alpha_n)/\Phi^*_{n+1}(\alpha_n)}$ and
\beq \label{en}
e_n=\sqrt{{\varpi_n(\alpha_n)\over\varpi_{n-1}(\alpha_{n-1})}{1 \over 1-|\lambda_{n-1}|^2}}.
\eeq

When $d\mu$ has an infinite support, there exists an infinite sequence of orthogonal
rational functions which generates an infinite sequence $(\lambda_n)$ in $\D$. If, on the
contrary, $d\mu$ is supported on a finite number $N+1$ of points, only the first $N+1$
orthogonal rational functions $\Phi_0,\dots,\Phi_N$ exist because $L^2(d\mu)$ is
$N+1$-dimensional. Nevertheless, there exists $\Phi_{N+1}\in\spn\{B_0,\dots,B_{N+1}\}$
and orthogonal to $\Phi_0,\dots,\Phi_N$, although it has $L^2(d\mu)$-norm equal to zero.
$\Phi_{N+1}$ satisfies a relation like (\ref{RR0}) with some coefficient $e_{N+1}\neq0$
and $\lambda_N\in\T$. In consequence, given a sequence $\bsalpha$ in $\OO$, we can
associate with any measure $d\mu$ a sequence $\bs\lambda=(\lambda_n)\in\frS$, which
terminates iff $d\mu$ is finitely supported, the number of points in the support being
equal to the length of $\bs\lambda$. $\bs\lambda$ will be called the
$\bsalpha$-parameters of $d\mu$. As follows from \cite[Theorem 8.1.4]{BGHN99}, this
establishes a surjective map
\beq \label{Ca-B}
\cS_\alpha \colon \kern0pt
\mathop{\frP_0 \xrightarrow{\kern0pt} \frS}
\limits_{\ds \hspace*{0pt} d\mu \to \bs\lambda}
\eeq
which is certainly bijective when $B_n\rightrightarrows0$. When $B_n$ does not diverge to
0, different probability measures can have the same $\bsalpha$-parameters.

Simultaneously, we can consider the S-function of $d\mu$ and the corresponding
$\bsalpha$-parameters $\bsgamma$. The key result is that $\bs\lambda=\bsgamma$, a fact
which is an immediate consequence of \cite[Corollary 6.5.2]{BGHN99} (to fit with the
notation there we must point out a misprint in formulas (6.27), (6.29) and (6.31), where
$z_n$ must be interchanged with $\overline z_n$; then, $\lambda_n=L_{n+1}$ and
$f_n=-\Gamma_n$). This result is the rational generalization of Geronimus' theorem for
the orthogonal polynomials on the unit circle (see \cite{Ge44,Ge54,Ge61}). Therefore, for
any sequence $\bsalpha=(\alpha_n)$ in $\OO$, we have the commutative diagram
$$
\begin{CD}
\frP_0              @>\cC_0>>          \frC_0   \\
@V\cS_\bsalpha VV   @VV\cB_0 V                  \\
\frS                @<<\cT_\bsalpha<   \frB
\end{CD}
$$
with $\cC_0$, $\cB_0$ homeomorphisms and $\cT_\bsalpha$, $\cS_\bsalpha$ continuous and
surjective. $\cS_\bsalpha$ and $\cT_\bsalpha$ are bijective when
\beq \label{BP}
B_n\rightrightarrows0 \kern5pt \Longleftrightarrow \kern5pt
\bca
\sum (1-|\alpha_n|) = \infty,
\smallskip \cr
\sum \frac{\im\,\alpha_n}{1+|\alpha_n|^2} = \infty,
\eca
\eeq
which means that not all the sequence $\bsalpha$ can approach to $\partial\OO$ very
quickly. If, on the contrary, $B_n$ does not diverge to 0, different S-functions can have
the same $\bsalpha$-parameters $\bsgamma$, and this occurs iff $\bsgamma$ does not determine
a unique probability measure. When this happens, we will say that we are in the indeterminate
case.

When $B_n\rightrightarrows0$, $\cS_\bsalpha$ and $\cT_\bsalpha$ are homeomorphisms, as can be
deduced from Theorem \ref{CONV-gamma} and the following important result.

\bt \label{CONV-gamma-inv}

Let $(\bsalpha^k)$ be a sequence of sequences in $\OO$, $\bsalpha$ a sequence in $\OO$
with associated Blaschke products $(B_n)$, $(f^k)$ a sequence in $\frB$ and $f\in\frB$.
If $\bsgamma^k$ are the $\bsalpha^k$-parameters of $f^k$ and $\bsgamma$ are the
$\bsalpha$-parameters of $f$, the condition $B_n\rightrightarrows0$ ensures that
$$
\bsalpha^k \to \bsalpha, \kern5pt \bsgamma^k \to \bsgamma
\kern5pt \Longrightarrow \kern5pt
f^k \rightrightarrows f.
$$

\et

\bpr
Let $\tilde f\in\Lim f^k$, i.e., $f^{k_j}\stackrel{j}{\rightrightarrows}\tilde f$ for
some subsequence $(f^{k_j})_j$. Then, $\tilde f \in \frB$ and we can consider its
$\bsalpha$-parameters $\tilde\bsgamma$. From Theorem \ref{CONV-gamma},
$\bsgamma^{k_j}\stackrel{j}{\to}\tilde\bsgamma$, thus $\tilde\bsgamma=\bsgamma$. Hence,
$B_n\rightrightarrows0$ ensures that $\tilde f=f$.
\epr

In Theorems \ref{CONV-gamma} and \ref{CONV-gamma-inv}, if $\bsgamma=(\gamma_k)_{k=0}^N$,
$N<\infty$, then the convergence condition $\bsalpha^k\to\bsalpha$ can be reduced to
$\alpha^k_n\stackrel{k}{\to}\alpha_n$ for $n \leq N+1$.

\bex \label{LEBESGUE2}

From example \ref{LEBESGUE}, $f(z;dm_\alpha)=\overline{\zeta_0(\alpha)}$ for each
$\alpha\in\OO$ and $f(z;\delta_\tau)=\overline{\zeta_0(\tau)}$ for any $\tau\in\partial\OO$.
So, no matter the sequence $\bsalpha$ in $\OO$,
$$
\cS_\bsalpha(dm_\alpha) = (\overline{\zeta_0(\alpha)},0,0,\dots),
\qquad
\cS_\bsalpha(\delta_\tau) = (\overline{\zeta_0(\tau)}).
$$
Consider a sequence $(\bsalpha^k)$ of sequences all in the same compact subset of $\OO$,
and a sequence $(d\mu^k)$ in $\frP_0$. If $\bsgamma^k=(\gamma^k_n)_n=\cS_{\bsalpha^k}(d\mu^k)$, then
$$
\ba{l}
\Lim d\mu^k \subset \{dm_\alpha : \alpha\in\OO\}
\kern5pt \Longrightarrow \kern5pt
\limsup|\gamma^k_0|<1, \kern5pt \gamma^k_n \stackrel{k}{\to} 0 \kern5pt \forall n\geq1,
\medskip \cr
\Lim d\mu^k \subset \{\delta_\tau : \tau\in\partial\OO\}
\kern5pt \Longrightarrow \kern5pt
|\gamma^k_0| \to 1.
\ea
$$
Besides, the above relations are ``iff" when the Blaschke product related to $\bsalpha$
diverges to 0. We will only prove the first right implication since the rest of them
follow analogous arguments.

Let us assume that $\Lim d\mu^k \subset \{dm_\alpha : \alpha\in\OO\}$ and let
$\bsgamma=(\gamma_n)\in\Lim\bsgamma^k$. Then, $\bsgamma^{k_j}\stackrel{j}{\to}\bsgamma$
for some subsequence $(\bsgamma^{k_j})_j$. Without loss of generality we can suppose
$d\mu^{k_j}\stackrel{*}{\underset{j}{\to}}dm_\alpha$, $\alpha\in\OO$, if necessary
restricting the subsequence. By a similar reason, we can assume that
$\bsalpha^{k_j}\stackrel{j}{\to}\bsalpha$ with $\bsalpha$ a sequence in the same compact
subset of $\OO$ than all the $\bsalpha^k$. From Theorems \ref{CONV-gral} and
\ref{CONV-gamma} we find that $\bsgamma^{k_j}\stackrel{j}{\to}\cS_\bsalpha(dm_\alpha)$,
thus $\gamma_0=\overline{\zeta_0(\alpha)}\in\D$ and $\gamma_n=0$ for $n\geq1$. Hence,
$\gamma^k_n\stackrel{k}{\to}0$ for $n\geq1$, and $\limsup|\gamma^k_0|<1$ because $\Lim\gamma^k_0$ is a
compact subset of $\D$.
{\huge$\hfill\centerdot$}

\eex

We finish this section showing that the relation between the $\alpha$-S-function and the
S-function of a measure leads to a connection between their $\bsalpha$-iterates and,
thus, between their $\bsalpha$-parameters. This connection is a simple consequence of the
following general results.

\bl \label{PROP-f}

Let $\bsalpha$ be a sequence in $\OO$, $f,\tilde f\in\frB$ and denote by $(f_n),(\tilde
f_n)$ and $(\gamma_n),(\tilde\gamma_n)$ the related $\bsalpha$-iterates and
$\bsalpha$-parameters respectively.
\be
\item
$\tilde f = \lambda f, \kern5pt \lambda\in\T
\kern5pt \Rightarrow \kern5pt
\tilde f_n = \lambda f_n, \kern5pt \tilde\gamma_n=\lambda\gamma_n, \kern5pt n\geq0.$
\item
$\ds \tilde f = \frac{f-w}{1-\overline wf}, \kern5pt w\in\D
\kern5pt \Rightarrow \kern5pt
\bca
\tilde f_0 = \frac{f_0-w}{1-\overline wf_0}, &
\tilde\gamma_0 = \frac{\gamma_0-w}{1-\overline w\gamma_0},
\smallskip \cr
\tilde f_n = \frac{1-w\overline\gamma_0}{1-\overline w\gamma_0}\,f_n, &
\tilde\gamma_n = \frac{1-w\overline\gamma_0}{1-\overline w\gamma_0}\,\gamma_n,
\kern12pt n\geq1.
\eca
$
\ee

\el

\bpr
The first item is trivial. For the second one, in view of the first result, it
suffices to prove it for $n=1$, which is just a matter of computation.
\epr

As a direct consequence of the previous result and Proposition \ref{f-fa} we find that
the sequences of $\bsalpha$-iterates and $\bsalpha$-parameters of S-functions and
$\alpha$-S-functions are proportional up to the first element of the sequence.

\bp \label{fn-fan}

Let $\alpha\in\OO$ and consider the S-function $f$ and the $\alpha$-S-function $f_\alpha$
of a measure $d\mu\in\frP_0$. If, for some sequence $\bsalpha$ in $\OO$,
$(f_n),(f_{\alpha,n})$ and $(\gamma_n),(\gamma_{\alpha,n})$ are the $\bsalpha$-iterates
and $\bsalpha$-parameters of $f,f_\alpha$ respectively,
$$
\ba{l} \ds
f_{\alpha,0} = -\frac{\zeta_0(\alpha)}{|\zeta_0(\alpha)|}
\frac{f_0-\overline{\zeta_0(\alpha)}}{1-\zeta_0(\alpha)f_0},
\qquad
\gamma_{\alpha,0} = -\frac{\zeta_0(\alpha)}{|\zeta_0(\alpha)|}
\frac{\gamma_0-\overline{\zeta_0(\alpha)}}{1-\zeta_0(\alpha)\gamma_0},
\medskip \cr \ds
f_{\alpha,n} = -\frac{\zeta_0(\alpha)}{|\zeta_0(\alpha)|}
\frac{1-\overline{\zeta_0(\alpha)\gamma_0}}{1-\zeta_0(\alpha)\gamma_0}\,f_n,
\kern10pt
\gamma_{\alpha,n} = -\frac{\zeta_0(\alpha)}{|\zeta_0(\alpha)|}
\frac{1-\overline{\zeta_0(\alpha)\gamma_0}}{1-\zeta_0(\alpha)\gamma_0}\,\gamma_n,
\kern10pt n\geq1.
\ea
$$

\ep

\section{Wall rational functions} \label{WR}

Consider a sequence $\bsalpha=(\alpha_n)$ in $\OO$ and a S-function $f$ with
$\bsalpha$-parameters $\bsgamma=(\gamma_n)$. The inverse relation between the
$\bsalpha$-iterates $(f_n)$ of $f$ can be written as  $f_{n-1}=M(\alpha_n,\gamma_{n-1})f_n$,
where
\beq \label{NP-inv}
M(\alpha,\gamma) f = \frac{\zeta_\alpha f+\gamma}{1+\overline\gamma\,\zeta_\alpha f} =
\gamma +\frac{(1-|\gamma|^2)\,\zeta_\alpha}{\overline\gamma\,\zeta_\alpha+\ds\frac{1}{f}}.
\eeq
The identity $f = M(\alpha_1,\gamma_0) M(\alpha_2,\gamma_1) \cdots M(\alpha_n,\gamma_{n-1}) f_n$
shows that
\beq \label{f-fn-CF}
\ba{l} \ds
f = \gamma_0 +
\frac{(1-|\gamma_0|^2)\,\zeta_1}{\overline\gamma_0\,\zeta_1} \underset{\ds+}{}
\frac{1}{\gamma_1} \underset{\ds+}{}
\frac{(1-|\gamma_1|^2)\,\zeta_2}{\overline\gamma_1\,\zeta_2} \underset{\ds+}{}
\cdots
\medskip \cr \ds \kern180pt
\cdots
\underset{\ds+}{} \frac{1}{\gamma_{n-1}} \underset{\ds+}{}
\frac{(1-|\gamma_{n-1}|^2)\,\zeta_n}{\overline\gamma_{n-1}\,\zeta_n} \underset{\ds+}{}
\frac{1}{f_n}.
\ea
\eeq
This provides a formal expansion of $f$ as an $\bsalpha$-dependent continued fraction
$$
f \sim \gamma_0 +
\frac{(1-|\gamma_0|^2)\,\zeta_1}{\overline\gamma_0\,\zeta_1} \underset{\ds+}{}
\frac{1}{\gamma_1} \underset{\ds+}{}
\frac{(1-|\gamma_1|^2)\,\zeta_2}{\overline\gamma_1\,\zeta_2} \underset{\ds+}{}
\cdots \underset{\ds+}{}
\frac{1}{\gamma_{n-1}} \underset{\ds+}{}
\frac{(1-|\gamma_{n-1}|^2)\,\zeta_n}{\overline\gamma_{n-1}\,\zeta_n} \underset{\ds+}{}
\cdots,
$$
which will be called the $\bsalpha$-continued fraction of $f$. Its $2n-2$ and $2n-1$
approximants will be denoted $f^{(n)}$ and $\tilde f^{(n)}$ respectively, i.e.,
$$
\ba{l} \ds
f^{(1)} = \gamma_0,
\quad
\tilde f^{(1)} = \gamma_0 +
\frac{(1-|\gamma_0|^2)\,\zeta_1}{\overline\gamma_0\,\zeta_1},
\quad
f^{(2)} = \gamma_0 +
\frac{(1-|\gamma_0|^2)\,\zeta_1}{\overline\gamma_0\,\zeta_1} \underset{\ds+}{}
\frac{1}{\gamma_1},
\medskip \cr \ds
\tilde f^{(2)} = \gamma_0 +
\frac{(1-|\gamma_0|^2)\,\zeta_1}{\overline\gamma_0\,\zeta_1} \underset{\ds+}{}
\frac{1}{\gamma_1} \underset{\ds+}{}
\frac{(1-|\gamma_1|^2)\,\zeta_2}{\overline\gamma_1\,\zeta_2},
\quad \dots
\ea
$$
Notice that susbstituting $f_n$ or $1/f_n$ by 0 in (\ref{f-fn-CF}) yields respectively
$f^{(n)}$ or $\tilde f^{(n)}$ instead of $f$. When $f$ has a finite sequence
$\bsgamma=(\gamma_n)_{n=0}^{N}$ of $\bsalpha$-parameters, the related continued fraction
is finite too because $\gamma_N\in\T$. In such a case, only the approximants
$f^{(1)},\dots,f^{(N+1)}$ and $\tilde f^{(1)},\dots,\tilde f^{(N)}$ exist, and
$f=f^{(N+1)}$ because the formal expansion as a continued fraction becomes an equality
since $f_N=\gamma_N$.

$M(\alpha,\gamma)$ transforms rational functions into rational functions, thus $f^{(n)}$
and $\tilde f^{(n)}$ are both rational functions. Moreover, if $\alpha\in\OO$ and
$\gamma\in\D$, $M(\alpha,\gamma)$ maps $\frB$ on $\frB^0$. Therefore,
$f^{(n)}\in\frB^0$ for all $n$, except for the case $N=0$ where
$f^{(1)}=f\in\frB\setminus\frB^0$. Two principal questions arise: What can we say about
the expression and properties of $f^{(n)}$ and $\tilde f^{(n)}$? Do they converge to $f$?
The first question will lead to the rational analogue of the Wall polynomials. As for the
convergence of $f^{(n)}$, an immediate answer emerges from the Nevalinna-Pick
homeomorphisms.

\bt \label{CONV-W}

Let $\bsalpha$ be a sequence in $\OO$ with related Blaschke products $(B_n)$, $f\in\frB$
and $f^{(n)}$ the $2n-2$ approximant of the associated $\bsalpha$-continued fraction. If
$\bsgamma=(\gamma_n)_{n=0}^N$ are the $\bsalpha$-parameters of $f$, the
$\bsalpha$-parameters of $f^{(n)}$ are
$\bsgamma^{(n)}=(\gamma_0,\dots,\gamma_{n-1},0,0,\dots)$ for $n<N+1$. When $N=\infty$,
the limit points of $(f^{(n)})$ are S-functions with $\bsalpha$-parameters $\bsgamma$. In
particular,
$$
B_n \rightrightarrows 0 \kern5pt \Longrightarrow \kern5pt f^{(n)} \rightrightarrows f.
$$

\et

\bpr
If $\alpha\in\OO$ and $\gamma\in\D$, $g=M(\alpha,\gamma)h\in\frB^0$ and satisfies
$g(\alpha)=\gamma$ for any $h\in\frB$. Therefore, if $n<N+1$,
$\gamma_0,\dots,\gamma_{n-1}\in\D$ and the relation
$f^{(n)}=M(\alpha_1,\gamma_0)M(\alpha_2,\gamma_1)\cdots M(\alpha_n,\gamma_{n-1})\,0$
shows that $\gamma_0=f^{(n)}(\alpha_1)$ and
$M(\alpha_1,\gamma_0)^{-1}f^{(n)}=M(\alpha_2,\gamma_1)\cdots M(\alpha_n,\gamma_{n-1})\,0$
is the first $\bsalpha$-iterate of $f^{(n)}$. By induction, the first $n$
$\bsalpha$-parameters of $f^{(n)}$ are $\gamma_0,\dots,\gamma_{n-1}$ and the $n$-th
$\bsalpha$-iterate of $f^{(n)}$ is 0. Hence, the rest of the $\bsalpha$-parameters of
$f^{(n)}$ are null. If $N=\infty$, then $\bsgamma^{(n)}\to\bsgamma$. Thus, the continuity
of $\cT_\bsalpha$ implies that the limit points of $(f^{(n)})$ must be $S$-functions with
$\bsalpha$-parameters $\bsgamma$. The condition $B_n\rightrightarrows0$ ensures that
$\cT_\bsalpha$ is a homeomorphism, hence $f^{(n)} \rightrightarrows f$.
\epr

Notice that $\bsgamma^{(n)}$ is the sequence of $\tilde\bsalpha$-parameters of $f^{(n)}$
whenever
$\tilde\bsalpha=(\alpha_1,\dots,\alpha_n,\tilde\alpha_{n+1},\tilde\alpha_{n+2},\dots)$,
no matter the choice of $\tilde\alpha_j\in\OO$ for $j>n$.

To analyze the nature of the approximants $f^{(n)}$ and $\tilde f^{(n)}$ we start writing
the relation $f_{n-1}=M(\alpha_n,\gamma_{n-1})f_n$ between the $\bsalpha$-iterates of $f$
in the way
$$
\begin{pmatrix} f_{n-1} \cr 1 \end{pmatrix} \doteq
\begin{pmatrix} \zeta_n & \gamma_{n-1} \cr \overline\gamma_{n-1}\zeta_n & 1 \end{pmatrix}
\begin{pmatrix} f_n \cr 1 \end{pmatrix},
$$
where the symbol $\doteq$ means equality up to a non vanishing scalar factor.
Therefore,
$$
\begin{pmatrix} f \cr 1 \end{pmatrix} \doteq
\begin{pmatrix} \zeta_1 & \gamma_0 \cr \overline\gamma_0\zeta_1 & 1 \end{pmatrix}
\begin{pmatrix} \zeta_2 & \gamma_1 \cr \overline\gamma_1\zeta_2 & 1 \end{pmatrix}
\cdots
\begin{pmatrix} \zeta_{n+1} & \gamma_n \cr \overline\gamma_n\zeta_{n+1} & 1 \end{pmatrix}
\begin{pmatrix} f_{n+1} \cr 1 \end{pmatrix}.
$$
It is evident that
$$
\begin{pmatrix} f \cr 1 \end{pmatrix} \doteq
\begin{pmatrix} \zeta_{n+1}\tilde S_n & R_n \cr \zeta_{n+1}\tilde R_n & S_n \end{pmatrix}
\begin{pmatrix} f_{n+1} \cr 1 \end{pmatrix},
$$
$R_n,S_n,\tilde R_n,\tilde S_n$ being linear combinations of the first $n+1$ Blaschke
products $B_0,B_1,\dots,B_n$ related to $\bsalpha$, with coefficients depending only on
the parameters $\gamma_0,\dots,\gamma_n$. Hence,
\beq \label{f-fn-tilde}
f = \frac{R_{n-1}+\tilde S_{n-1}\,\zeta_nf_n}{S_{n-1}+\tilde R_{n-1}\,\zeta_nf_n},
\eeq
which is a compact way of writing (\ref{f-fn-CF}). As we mentioned before,
substituting $f_{n+1}$ or $1/f_{n+1}$ by 0 in (\ref{f-fn-CF}), i.e. in (\ref{f-fn-tilde}),
we get respectively $f^{(n)}$ or $\tilde f^{(n)}$ instead of $f$. Thus,
$$
f^{(n)}=\frac{R_{n-1}}{S_{n-1}}, \qquad \tilde f^{(n)}=\frac{\tilde S_{n-1}}{\tilde R_{n-1}}.
$$
Besides, $\tilde R_n$ and $\tilde S_n$ can be expressed in terms of $R_n$ and $S_n$.
From the equality
$$
\begin{pmatrix} \zeta_{n+1}\tilde S_n & R_n \cr \zeta_{n+1}\tilde R_n & S_n \end{pmatrix} =
\begin{pmatrix} \zeta_n\tilde S_{n-1} & R_{n-1} \cr \zeta_n\tilde R_{n-1} & S_{n-1} \end{pmatrix}
\begin{pmatrix} \zeta_{n+1} & \gamma_n \cr \overline\gamma_n\zeta_{n+1} & 1 \end{pmatrix},
$$
we obtain
\beq \label{RR-W}
\begin{pmatrix} \tilde S_n & R_n \cr \tilde R_n & S_n \end{pmatrix} =
\begin{pmatrix} \tilde S_{n-1} & R_{n-1} \cr \tilde R_{n-1} & S_{n-1} \end{pmatrix}
\begin{pmatrix} \zeta_n & \gamma_n\zeta_n \cr \overline\gamma_n & 1 \end{pmatrix},
\eeq
which, together with the initial condition
\beq \label{CI-W}
\begin{pmatrix} \tilde S_0 & R_0 \cr \tilde R_0 & S_0 \end{pmatrix} =
\begin{pmatrix} 1 & \gamma_0 \cr \overline\gamma_0 & 1 \end{pmatrix},
\eeq
permits us to prove by induction that $\tilde R_n=R^*_n=B_nR_{n*}$ and $\tilde S_n=S^*_n=B_nS_{n*}$.
In consequence,
\beq \label{f-fn}
f = \frac{R_{n-1}+S^*_{n-1}\,\zeta_nf_n}{S_{n-1}+R^*_{n-1}\,\zeta_nf_n},
\qquad
f^{(n)}=\frac{R_{n-1}}{S_{n-1}},
\qquad
\tilde f^{(n)}=\frac{S^*_{n-1}}{R^*_{n-1}}.
\eeq
where $R_n$ and $S_n$ are recursively defined by
\beq \label{RR-W2}
\bca
R_0=\gamma_0,
\smallskip \cr
S_0=1,
\eca
\qquad
\bca
R_n=R_{n-1}+\gamma_n\zeta_nS^*_{n-1},
\smallskip \cr
S_n=S_{n-1}+\gamma_n\zeta_nR^*_{n-1},
\eca
\quad n\geq1.
\eeq
We will call $(R_n)$ and $(S_n)$ the Wall rational functions associated with $f$ and
$\bsalpha$. The number of Wall rational functions coincides with the number of
$\bsalpha$-parameters of $f$. The recurrence for $R^*_n$ and $S^*_n$
\beq \label{RR-W*}
\bca
R^*_0=\overline\gamma_0,
\smallskip \cr
S^*_0=1,
\eca
\qquad
\bca
R^*_n=\zeta_nR^*_{n-1}+\overline\gamma_nS_{n-1},
\smallskip \cr
S^*_n=\zeta_nS^*_{n-1}+\overline\gamma_nR_{n-1},
\eca
\quad n\geq1,
\eeq
shows that $S^*_n$ is a monic element of $\spn\{B_0,\dots,B_n\}$, i.e., the coefficient
of $B_n$ in the expansion of $S^*_n$ as a linear combination of $B_0,\dots,B_n$ is 1.
We can also get from (\ref{RR-W}) the inverse recurrence
$$
\bca
R_{n-1}=(1-|\gamma_n|^2)^{-1}(R_n-\gamma_nS^*_n),
\smallskip \cr
S_{n-1}=(1-|\gamma_n|^2)^{-1}(S_n-\gamma_nR^*_n),
\eca
\quad 1 \leq n < N.
$$

The expression of the approximants $f^{(n)}$ and $\tilde f^{(n)}$ in terms of the Wall
rational functions implies that $\tilde f^{(n)}=1/f^{(n)}_*$. Thus, we obtain the
following result for the convergence of $\tilde f^{(n)}$ as a direct consequence of
Theorem \ref{CONV-W}.

\bt \label{CONV-W*}

Let $\bsalpha$ be a sequence in $\OO$ with related Blaschke products $(B_n)$, $f\in\frB$
and $\tilde f^{(n)}$ the $2n-1$ approximant of the associated $\bsalpha$-continued
fraction. If $f$ has an infinite sequence of $\bsalpha$-parameters,
$$
B_n \rightrightarrows 0
\kern5pt \Longrightarrow \kern5pt
\tilde f^{(n)} \rightrightarrows 1/f_*
\text{ in } \OO^e\setminus\{z\in\OO^e : f(\hat z)=0\}.
$$
\et

The above result has a natural interpretation in view of the integral representation of
any S-function $f$. Let $d\mu$ be the probability measure such that $f(z)=f(z;d\mu)$. The
expression for the corresponding C-function $F(z)=F(z;d\mu)$ defines an analytic
function, not only for $z\in\OO$, but also for $z\in\OO^e$. This permits us to extend the
definition of $f$ to the points of $\OO^e$ throughout
$f=\frac{1}{\zeta_0}\frac{F-1}{F+1}$. It is direct to see that $F=-F_*$, thus $f=1/f_*$.
This means that the extended $f$ is analytic at $z\in\OO^e$ iff $f$ has not a zero at
$\hat z\in\OO$. With such an extended $f$, Theorems \ref{CONV-W} and \ref{CONV-W*} can be
combined to saying that, when $B_n\rightrightarrows0$, then $f^{(n)} \rightrightarrows f$
in $\OO$, while $\tilde f^{(n)} \rightrightarrows f$ in $\{z\in\OO^e : f \text{ is
analytic at } z\}$.

We can state some general properties of the Wall rational functions.

\bp \label{PROP-W}

Let $\bsalpha$ be a sequence in $\OO$ with Blaschke products $(B_n)$, $(R_n)_{n=0}^N$,
$(S_n)_{n=0}^N$ the Wall rational functions associated with $f\in\frB$, $(f_n)_{n=0}^N$
the corresponding $\bsalpha$-iterates and $(\gamma_n)_{n=0}^N$ the related
$\bsalpha$-parameters. Let us denote $\Upsilon_n=\prod_{k=0}^n(1-|\gamma_k|^2)$.

\be

\item
$S_n(z) \, \overline{S_n(w)} - R^*_n(z) \, \overline{R^*_n(w)} =$

\hskip50pt $=(1-|\gamma_n|^2) (S_{n-1}(z) \, \overline{S_{n-1}(w)} - \zeta_n(z) \,
\overline{\zeta_n(w)} \, R^*_{n-1}(z) \, \overline{R^*_{n-1}(w)})$.

\item
$|S_n|^2-|R^*_n|^2 \geq \Upsilon_n$ in $\overline\OO$, the inequality being an equality
in $\partial\OO$. Furthermore, $\Upsilon_n^{-1}(|S_n|^2-|R^*_n|^2)$ is a
non-decreasing sequence in $\overline\OO$.

\item
$S_nS^*_n-R_nR^*_n = \Upsilon_nB_n$.

\item
$S_n$ does not vanish in $\overline\OO$ and has no zeros in common with $R_n$, neither
with $R^*_n$.

\item
$\ds \frac{R_n}{S_n}, \frac{R^*_n}{S_n}, \frac{\Upsilon_n}{S^2_n} \in \frB$. Moreover,
$\ds \left|\frac{R_n}{S_n}\right|, \left|\frac{R^*_n}{S_n}\right| <1$ in $\overline\OO$ if $n<N$
and \break
\vskip-20pt
$\ds \frac{R_N}{S_N}=f, \frac{R^*_N}{S_N}=\overline\gamma_N$ if $N<\infty$.

\item
$\ds f-\frac{R_n}{S_n} = \frac{\Upsilon_n}{S^2_n}
\frac{B_{n+1}f_{n+1}}{1+\frac{R^*_n}{S_n}\,\zeta_{n+1}f_{n+1}}$.

\item
$\ds \left|f-\frac{R_n}{S_n}\right| \leq (1+|\zeta_{n+1}|)|B_{n+1}|$.

\item
$\bca
S_n+R^*_n\zeta_{n+1}f_{n+1} = \prod_{k=0}^n(1+\overline\gamma_k\zeta_{k+1}f_{k+1}),
\smallskip \cr
R_n+S^*_n\zeta_{n+1}f_{n+1} = (\gamma_0+\zeta_1f_1)\prod_{k=1}^n(1+\overline\gamma_k\zeta_{k+1}f_{k+1}).
\eca$

\item
$\bca
S_nf-R_n = B_{n+1}f_{n+1}\prod_{k=0}^n(1-\overline\gamma_kf_k),
\smallskip \cr
S^*_n-R^*_nf = B_n\prod_{k=0}^n(1-\overline\gamma_kf_k). \eca$

\ee

\ep

\bpr
Property 1 is a direct consequence of (\ref{RR-W2}) and (\ref{RR-W*}).
Evaluating it for
$w=z$ gives $|S_n|^2-|R^*_n|^2 = (1-|\gamma_n|^2)
(|S_{n-1}|^2-|\zeta_n|^2|R^*_{n-1}|^2)$. Thus, $|S_n|^2-|R^*_n|^2 \geq (1-|\gamma_n|^2)
(|S_{n-1}|^2-|R^*_{n-1}|^2)$ in $\overline\OO$, with an equality in $\partial\OO$. This
proves 2 since $|S_0|^2-|R^*_0|^2=1-|\gamma_0|^2$.

Taking determinants in (\ref{RR-W}) and (\ref{CI-W}) we get
$S_0S^*_0-R_0R^*_0=1-|\gamma_0|^2$ and
$S_nS^*_n-R_nR^*_n=\zeta_n(1-|\gamma_n|^2)(S_{n-1}S^*_{n-1}-R_{n-1}R^*_{n-1})$, which
proves 3. The last equality follows also from Property 1 for $w=\hat z$.

Property 2 implies that, for $n<N$, $|S_n|^2\geq\Upsilon_n>0$ and $|S_n|>|R^*_n|$ in
$\overline\OO$. Besides, if $N<\infty$, using (\ref{RR-W2}) we get
$|S_N|\geq|S_{N-1}|-|R^*_{N-1}|>0$ in $\overline\OO$. Hence, $S_n$ does not vanish in
$\overline\OO$ for any $n$. Then, Property 3 shows that $S_n$ has no common zeros with
$R_nR^*_n$ because $B_n$ only vanishes at $\alpha_1,\dots,\alpha_n\in\OO$.

Since $S_n$ does not vanish in $\overline\OO$, $R_n/S_n$, $R^*_n/S_n$ and
$\Upsilon_n/S^2_n$ are analytic in a neighbourhood of $\overline\OO$. From 2,
$|\Upsilon_n/S^2_n|\leq1$ and $|R^*_n/S_n|<1$ for $n<N$ in $\overline\OO$, while
(\ref{RR-W2}) and (\ref{RR-W*}) give $R^*_N/S_N=\overline\gamma_N\in\T$ when $N<\infty$.
Also, we know that $R_n/S_n=f^{(n)}\in\frB^0$ if $n<N$ and $R_N/S_N=f^{(N+1)}=f\in\frB$
for $N<\infty$. Moreover, $|R_n/S_n|=|R^*_n/S_n|<1$ in $\partial\OO$ for $n<N$.

Property 6 is obtained using (\ref{f-fn}) to express $f-f^{(n+1)}$ and simplifying the
result with Property 3.

To prove 7, write 6 in the way
$$
f-\frac{R_n}{S_n} = \frac{\Upsilon_n}{S^2_n}
\frac{B_{n+1}f_{n+1}}{1-\left(\frac{R^*_n}{S_n}\,\zeta_{n+1}f_{n+1}\right)^2}
\left(1-{\ts\frac{R^*_n}{S_n}}\,\zeta_{n+1}f_{n+1}\right).
$$
Then, use 2 and 5 to get $|\Upsilon_n/S^2_n| \leq 1-|R^*_n/S_n|^2 \leq
1-\left|\frac{R^*_n}{S_n}\,\zeta_{n+1}f_{n+1}\right|^2$ and
$\left|1-\frac{R^*_n}{S_n}\,\zeta_{n+1}f_{n+1}\right| \leq 1+|\zeta_{n+1}|$ in $\OO$.

From (\ref{NP-ITERATES}), (\ref{RR-W2}), (\ref{RR-W*}) and the help of (\ref{ID}) we arrive
at the identities
$S_n+R^*_n\zeta_{n+1}f_{n+1}=(1+\overline\gamma_n\zeta_{n+1}f_{n+1})(S_{n-1}+R^*_{n-1}\zeta_nf_n)$
and
$R_n+S^*_n\zeta_{n+1}f_{n+1}=(1+\overline\gamma_n\zeta_{n+1}f_{n+1})(R_{n-1}+S^*_{n-1}\zeta_nf_n)$.
This, together with the initial conditions given in (\ref{RR-W2}) and (\ref{RR-W*}),
proves 8 by induction.

Using Properties 3 and 8 in (\ref{f-fn}), and taking into account identity (\ref{ID}), we
get 9.
\epr

Property 7 of the above proposition measures the rate of the convergence
$R_n/S_n\rightrightarrows f$ when $B_n\rightrightarrows0$. In the polynomial situation it
yields $|f-R_n/S_n|\leq(1+|z|)|z|^{n+1}$, which improves the usual bounds given in the
literature for this case.

The recurrence for the Wall rational functions permits us to identify certain iterates of
the S-function $R^*_n/S_n$.

\bp \label{R*n/Sn}

Let $(R_n)$, $(S_n)$ be the Wall rational functions associated with a sequence
$\bsalpha=(\alpha_n)$ in $\OO$ and an S-function $f$ with $\bsalpha$-parameters
$\bsgamma=(\gamma_n)$. If
$\tilde\bsalpha=(\alpha_n,\alpha_{n-1},\dots,\alpha_1,\alpha_0,\alpha_0,\alpha_0,\dots)$,
then the $\tilde\bsalpha$-iterates and $\tilde\bsalpha$-parameters of $R^*_n/S_n$ are
respectively
$$
(R_n/S_n,R_{n-1}/S_{n-1},\dots,R_0/S_0,0,0,\dots), \qquad
(\overline\gamma_n,\overline\gamma_{n-1},\dots,\overline\gamma_0,0,0,\dots).
$$

\ep

\bpr
It is simply a consequence of the identity
\beq \label{R*S-schur}
\frac{R^*_n}{S_n} =
\frac{\zeta_nR^*_{n-1}/S_{n-1}+\overline\gamma_n}{1+\gamma_n\zeta_nR^*_{n-1}/S_{n-1}},
\eeq
which is obtained directly from (\ref{RR-W2}) and (\ref{RR-W*}).
\epr

Proposition \ref{R*n/Sn} also works with
$\tilde\bsalpha=(\alpha_n,\alpha_{n-1},\dots,\alpha_1,\tilde\alpha_{n+1},\tilde\alpha_{n+2},\dots)$,
where $\tilde\alpha_j$ are arbitrary points of $\OO$ for $j>n$. As an immediate
consequence of the previous results and Theorems \ref{CONV-gamma}, \ref{CONV-gamma-inv},
we have that, for any sequences $\bsalpha^k$, $\bsalpha$ in $\OO$ and any S-functions
$f^k$, $f$,
$$
\bsalpha^k \to \bsalpha, \kern5pt f^k \rightrightarrows f
\kern5pt \Longrightarrow \kern5pt
\frac{R^k_n}{S^k_n} \stackrel{k}{\rightrightarrows} \frac{R_n}{S_n}, \kern5pt
\frac{R^{k*}_n}{S^k_n} \stackrel{k}{\rightrightarrows} \frac{R^*_n}{S_n} \kern5pt
\forall n,
$$
where $(R^k_n)_n$, $(S^k_n)_n$ are the Wall rational functions associated with $f^k$,
$\bsalpha^k$ and $(R_n)$, $(S_n)$ are the Wall rational functions associated with $f$,
$\bsalpha$. Indeed, a stronger result can be obtained.

\bp \label{CONV-RS}

Let $\bsalpha^k$ be a sequence of sequences in $\OO$, $\bsalpha$ a sequence in $\OO$,
$(f^k)$ a sequence in $\frB$ and $f\in\frB$. If $(R^k_n)_n$, $(S^k_n)_n$ are the Wall
rational functions associated with $f^k$, $\bsalpha^k$ and $(R_n)$, $(S_n)$ are the Wall
rational functions associated with $f$, $\bsalpha$, then, for all $n$,
$$
\bsalpha^k \to \bsalpha, \kern5pt f^k \rightrightarrows f
\kern5pt \Longrightarrow \kern5pt
\bca
R^k_n \stackrel{k}{\rightrightarrows} R_n, \kern5pt
R^{k*}_n \stackrel{k}{\rightrightarrows} R^*_n,
\smallskip \cr
S^k_n \stackrel{k}{\rightrightarrows} S_n, \kern5pt
S^{k*}_n \stackrel{k}{\rightrightarrows} S^*_n,
\eca
\text{in } \C\setminus\{\hat\alpha_1,\dots,\hat\alpha_n\}.
$$

\ep

\bpr
In view of Theorem \ref{CONV-gamma}, $\bsalpha^k \to \bsalpha$ and $f^k \rightrightarrows
f$ imply $\bsgamma^k\to\bsgamma$, where $\bsgamma^k$ are the $\bsalpha^k$-parameters of
$d\mu^k$ and $\bsgamma$ are the $\bsalpha$-parameters of $d\mu$. Then, the proof follows
by induction using (\ref{RR-W2}), (\ref{RR-W*}) and taking into account that
$\zeta_{\alpha^k_n} \stackrel{k}{\rightrightarrows} \zeta_{\alpha_n}$ in
$\C\setminus\{\hat\alpha_n\}$.
\epr

Given a sequence $\bsalpha$ in $\OO$, the Wall rational functions $(R_n)$, $(S_n)$
associated with an S-function $f$ are related to the orthogonal rational functions
$(\Phi_n)$ corresponding to the measure $d\mu$ such that $f(z)=f(z;d\mu)$. The relation
also involves the so called second kind rational functions $(\Psi_n)$, defined by
$$
\ba{l}
\Psi_0(z)=1,
\smallskip \cr
\Psi_n(z)=\int D(t,z) (\Phi_n(t)-\Phi_n(z)) \, d\mu(t), \qquad n\geq1.
\ea
$$
$(\Psi_n)$ are orthogonal rational functions associated with the same sequence
$\bsalpha$, but with respect to a measure with $\bsalpha$-parameters opposed to those
ones of $d\mu$ (see \cite[Theorems 4.2.4 and 6.2.5]{BGHN99}). Therefore,
$(\Psi_n,-\Psi^*_n)$ satisfy the same recurrence (\ref{RR0}) as $(\Phi_n,\Phi^*_n)$, but
with a different initial condition.

\bp \label{W-ORF}

Let $\bsalpha$ be a sequence in $\OO$, $(R_n)_{n=0}^N$, $(S_n)_{n=0}^N$ the Wall rational
functions related to $f\in\frB$ and $(\Phi_n)_{n=0}^N$, $(\Psi_n)_{n=0}^N$ the orthogonal
and second kind rational functions for the measure $d\mu\in\frP_0$ such that
$f(z)=f(z;d\mu)$. Denoting
$\kappa_n=(\Upsilon_{n-1}\varpi_0(\alpha_0)/\varpi_n(\alpha_n))^{1/2}$, we have for
$n<N+1$,
$$
\ba{l}
\bca
R_{n-1} = \frac{\kappa_n}{2} \frac{\varpi_n}{\varpi^*_0} (\Psi^*_n-\Phi^*_n),
\smallskip \cr
R^*_{n-1} = \frac{\kappa_n}{2z_n} \frac{\varpi_n}{\varpi_0} (\Psi_n-\Phi_n),
\eca
\kern39pt
\bca
S_{n-1} = \frac{\kappa_n}{2} \frac{\varpi_n}{\varpi_0} (\Psi^*_n+\Phi^*_n),
\smallskip \cr
S^*_{n-1} = \frac{\kappa_n}{2z_n} \frac{\varpi_n}{\varpi^*_0} (\Psi_n+\Phi_n),
\eca
\bigskip \cr
\bca
\Phi_n = \frac{z_n}{\kappa_n} \frac{\varpi_0}{\varpi_n} (\zeta_0S^*_{n-1}-R^*_{n-1}),
\smallskip \cr
\Phi^*_n = \frac{1}{\kappa_n} \frac{\varpi_0}{\varpi_n} (S_{n-1}-\zeta_0R_{n-1}),
\eca
\qquad
\bca
\Psi_n = \frac{z_n}{\kappa_n} \frac{\varpi_0}{\varpi_n} (\zeta_0S^*_{n-1}+R^*_{n-1}),
\smallskip \cr
\Psi^*_n = \frac{1}{\kappa_n} \frac{\varpi_0}{\varpi_n} (S_{n-1}+\zeta_0R_{n-1}).
\eca
\ea
$$

\ep

\bpr
Let us denote $\hat\Phi_n=\overline z_n\Phi_n$, $\hat\Psi_n=\overline z_n\Psi_n$ and
$$
\W_n =
\begin{pmatrix}
S^*_n & -R^*_n
\cr
-R_n & S_n \end{pmatrix},
\qquad
\F_n =
\begin{pmatrix}
\hat\Phi_n & \hat\Psi_n
\cr
\Phi^*_n & -\Psi^*_n
\end{pmatrix}.
$$
With this notation, the recurrences for the Wall rational functions and the orthogonal
and second kind rational functions read as
$$
\W_n=T_n\W_{n-1}, \qquad \F_n=e_n\frac{\varpi_{n-1}}{\varpi_n}\,T_{n-1}\F_{n-1},
$$
with $T_n$ and $e_n$ given in (\ref{RR}) and (\ref{en}) respectively. Hence,
$$
\ba{l}
\W_n = T_n \cdots T_1 \W_0 =
T_n \cdots T_1 T_0 \begin{pmatrix} \zeta_0 & 0 \cr 0 & 1 \end{pmatrix}^{\kern-3pt -1},
\medskip \cr \ds
\F_n = \frac{1}{\kappa_n} \frac{\varpi_0}{\varpi_n} \,
T_{n-1} \cdots T_1 T_0 \F_0 =
T_n \cdots T_1 T_0 \begin{pmatrix} 1 & 1 \cr 1 & -1 \end{pmatrix}.
\ea
$$
So,
$$
\F_n = \frac{1}{\kappa_n} \frac{\varpi_0}{\varpi_n} \, \W_{n-1}
\begin{pmatrix} \zeta_0 & \zeta_0 \cr 1 & -1 \end{pmatrix},
$$
which gives the desired relations.
\epr

The above result allows us to identify the measure and C-function corresponding to the
S-function $f^{(n)}=R_{n-1}/S_{n-1}$. The measures $dm_\alpha$, $\alpha\in\OO$, given
in Example \ref{LEBESGUE} play an important role in such an identification.

\bp \label{ALL-BERNSTEIN}

Let $\bsalpha$ be a sequence in $\OO$, $(\Phi_n)_{n=0}^N$, $(\Psi_n)_{n=0}^N$ the
orthogonal and second kind rational functions of $d\mu\in\frP_0$, $(R_n)_{n=0}^N$,
$(S_n)_{n=0}^N$ the Wall rational functions related to $f(z;d\mu)$ and
$\bsgamma=(\gamma_n)_{n=0}^N$ its $\bsalpha$-parameters. Then, for $n<N+1$, we have the
correspondences
$$
\mathop{
\bcd
\frP_0   @>\cC_0>>   \frC_0   @>\cB_0>>   \frB   @>\cT_\bsalpha>>   \frS
\ecd
}
\limits_{\ds \kern80pt
\frac{\,dm_{\alpha_n}}{\,|\Phi_n|^2} \xrightarrow{\kern17pt}
\frac{\Psi^*_n}{\Phi^*_n} \xrightarrow{\kern18pt}
\frac{R_{n-1}}{S_{n-1}} \xrightarrow{\kern5pt}
(\gamma_0,\dots,\gamma_{n-1},0,0,\dots)
}
$$

\ep

\bpr
From Theorem \ref{CONV-W} we know that $(\gamma_0,\dots,\gamma_{n-1},0,0,\dots)$ are the
$\bsalpha$-parameters of $f^{(n)}=R_{n-1}/S_{n-1}$. Besides, Theorem \ref{W-ORF} gives
$$
\frac{\Psi^*_n}{\Phi^*_n} =
\frac{1+\zeta_0R_{n-1}/S_{n-1}}{1-\zeta_0R_{n-1}/S_{n-1}},
$$
which shows that $\Psi^*_n/\Phi^*_n$ is a C-function and
$R_{n-1}/S_{n-1}=\cB_0(\Psi^*_n/\Phi^*_n)$. Finally, the fact that
$dm_{\alpha_n}/|\Phi_n|^2$ is a probability measure with C-function $\Psi^*_n/\Phi^*_n$
was proven in \cite[Theorem 4.2.6]{BGHN99}.
\epr

Given a sequence $\bsalpha$ in $\OO$ and a measure $d\mu\in\frP_0$ with orthogonal
rational functions $(\Phi_n)$, we will denote
\beq \label{mun}
d\mu^{(n)}=\frac{\,dm_{\alpha_n}}{\,|\Phi_n|^2},
\eeq
so that, according to the previous notation, if $f$ is the S-function of $d\mu$,
$\bsgamma=(\gamma_n)$ its $\bsalpha$-parameters and $(\Psi_n)$ the second kind rational functions,
$$
\ba{c}
\ds F(z;d\mu^{(n)})=\frac{\Psi^*_n(z)}{\Phi^*_n(z)},
\qquad
f(z;d\mu^{(n)})=\frac{R_{n-1}(z)}{S_{n-1}(z)}=f^{(n)}(z),
\medskip \cr
\cS_\bsalpha(d\mu^{(n)})=(\gamma_0,\dots,\gamma_{n-1},0,0,\dots)=\bsgamma^{(n)}.
\ea
$$
Notice that, if
$\tilde\bsalpha=(\alpha_1,\dots,\alpha_n,\tilde\alpha_{n+1},\tilde\alpha_{n+2},\dots)$
with $\tilde\alpha_j$ arbitrary points of $\OO$ for $j>n$, then
$\cS_{\tilde\bsalpha}(d\mu^{(n)})=\bsgamma^{(n)}$ too. Using recurrence (\ref{RR0}) we see
that the orthogonal rational functions associated with $d\mu^{(n)}$ and $\tilde\bsalpha$
are $(\Phi_0,\dots,\Phi_n,\tilde\Phi_{n+1},\tilde\Phi_{n+2},\dots)$ where, using the
tilde to refer to the the elements related to $\tilde\bsalpha$,

$$
\tilde\Phi_j = \sqrt{\frac{\tilde\varpi_j(\tilde\alpha_j)}{\varpi_n(\alpha_n)}}
\frac{\varpi^*_n}{\tilde\varpi_j} \frac{\tilde B_{j-1}}{B_n} \, \Phi_n, \qquad j>n.
$$

Remember that $R_n/S_n$ and $S^*_n/R^*_n$ are respectively the $2n$ and $2n+1$
approximants of the $\bsalpha$-continued fraction for $f(z)=f(z;d\mu)$. In \cite[Section
4.4]{BGHN99} it is shown that, analogously, $\Psi^*_n/\Phi^*_n$ and $-\Psi_n/\Phi_n$ are
respectively the $2n$ and $2n+1$ approximants of an $\bsalpha$-dependent continued
fraction expansion of $F(z;d\mu)$ given by
$$
\ba{l} \ds
1 - \frac{2}{1} \underset{\ds+}{} \frac{-1}{\gamma_0\,\zeta_0} \underset{\ds+}{}
\frac{(1-|\gamma_0|^2)\,\zeta_0}{\overline\gamma_0} \underset{\ds+}{}
\frac{1}{\gamma_1\,\zeta_1} \underset{\ds+}{}
\frac{(1-|\gamma_1|^2)\,\zeta_1}{\overline\gamma_1} \underset{\ds+}{}
\cdots
\medskip \cr \ds \kern210pt
\cdots
\underset{\ds+}{} \frac{1}{\gamma_n\,\zeta_n} \underset{\ds+}{}
\frac{(1-|\gamma_n|^2)\,\zeta_n}{\overline\gamma_n} \underset{\ds+}{}
\cdots,
\ea
$$
so that the even and odd approximants converge to $F(z;d\mu)$ in $\OO$ and $\OO^e$
respectively when the Blaschke product related to $\bsalpha$ diverges to 0. Under this
condition we also have $d\mu^{(n)}\stackrel{*}{\to}d\mu$.

\section{Rational Khrushchev's formula} \label{KF}

The preceding results allow us to obtain a rational analogue of Khrushchev's formula for
the orthogonal polynomials on the unit circle (see \cite[Theorems 2 and 3]{Kh01}). We
will state first a weak version of it. In what follows, the $\bsalpha$-iterates of $d\mu$
means the $\bsalpha$-iterates of $f(z;d\mu)$.

\bt \label{K1}

Let $\bsalpha=(\alpha_n)$ be a sequence in $\OO$, $d\mu\in\frP_0$ and $(\Phi_n)$ the
related orthogonal rational functions. If $b_n=\overline z_n\Phi_n/\Phi^*_n$ and $(f_n)$
are the $\bsalpha$-iterates of $d\mu$, then
$$
\frac{|\Phi_n(t)|^2}{D_R(t,\alpha_n)}\,\mu'(t) =
\frac{1-|f_n(t)|^2}{|1-\zeta_n(t)b_n(t)f_n(t)|^2},
\qquad \text{ a.e. } t\in\partial\OO.
$$

\et

\bpr
From (\ref{REF}) we find that
$$
\mu' = \re\left(\frac{1+\zeta_0f}{1-\zeta_0f}\right) = \frac{1-|f|^2}{|1-\zeta_0f|^2},
\qquad \text{a.e. on } \partial\OO.
$$
(\ref{f-fn}) and Property 2 of Proposition \ref{PROP-W} yield
$$
1-|f|^2 = \Upsilon_{n-1} \frac{1-|f_n|^2}{|S_{n-1}+R^*_{n-1}\zeta_nf_n|^2},
\qquad \text{a.e. on } \partial\OO.
$$
Using again (\ref{f-fn}), together with the relations of Proposition \ref{W-ORF}, we obtain
$$
1-\zeta_0f = \kappa_n \frac{\varpi_n}{\varpi_0}
\frac{\Phi^*_n-\overline z_n\zeta_n\Phi_nf_n}{S_{n-1}+R^*_{n-1}\zeta_nf_n},
$$
hence
$$
|1-\zeta_0f|^2 =
\Upsilon_{n-1} \frac{\varpi_0(\alpha_0)}{\varpi_n(\alpha_n)} \left|\frac{\varpi_n}{\varpi_0}\right|^2
\frac{|\Phi_n|^2|1-\zeta_nb_nf_n|^2}{|S_{n-1}+R^*_{n-1}\zeta_nf_n|^2},
\qquad \text{a.e. on } \partial\OO.
$$
Combining the previous equalities and taking into account (\ref{RED}) we get the result.
\epr

Notice that the equality of Theorem \ref{K1} is trivial when $d\mu$ is finitely supported
because then $f_n$ is a finite Blaschke product.

The functions $b_n=\overline z_n\Phi_n/\Phi^*_n$ are finite Blaschke products because the
zeros of $\Phi_n$ lie on $\OO$. Concerning their iterates, we have the following result.

\bp \label{bn}

Let $b_n=\overline z_n\Phi_n/\Phi^*_n$, where $(\Phi_n)$ are the orthogonal rational
functions associated with a sequence $\bsalpha=(\alpha_n)$ in $\OO$ and a measure
$d\mu\in\frP_0$ with $\bsalpha$-parameters $\bsgamma=(\gamma_n)$. If
$\tilde\bsalpha=(\alpha_{n-1},\alpha_{n-2},\dots,\alpha_1,\alpha_0,\alpha_0,\alpha_0,\dots)$,
then, the $\tilde\bsalpha$-iterates and $\tilde\bsalpha$-parameters of $b_n$ are
respectively
$$
(b_n,b_{n-1},\dots,b_0), \qquad
(-\overline\gamma_{n-1},-\overline\gamma_{n-2},\dots,-\overline\gamma_0,1).
$$

\ep

\bpr
It follows immediately from the identity
$$
b_n = \frac{\zeta_{n-1}b_{n-1}-\overline\gamma_{n-1}}{1-\gamma_{n-1}\zeta_{n-1}b_{n-1}}
$$
obtained from recurrence (\ref{RR0}).
\epr

The above result also holds if
$\tilde\bsalpha=(\alpha_{n-1},\alpha_{n-2},\dots,\alpha_0,\tilde\alpha_{n+1},\tilde\alpha_{n+2},\dots)$,
where $\tilde\alpha_j$ are arbitrary points of $\OO$ for $j>n$. Following Khrushchev's
terminology, we will call $(b_n)$ the sequence of inverse $\bsalpha$-iterates of
$f(z;d\mu)$ or, equivalently, of $d\mu$. From the above proposition and Theorems
\ref{CONV-gamma}, \ref{CONV-gamma-inv}, we easily get a convergence property for the
inverse $\bsalpha^k$-iterates of $d\mu^k$ when $\bsalpha^k\to\bsalpha$ and $d\mu^k
\stackrel{*}{\to} d\mu$. Moreover, using the relations of Proposition \ref{W-ORF}, we
obtain from Proposition \ref{CONV-RS} a similar convergence property for the orthogonal
rational functions of $d\mu^k$. We summarize all these results.

\bp \label{CONV-ORF}

Let $(\bsalpha^k)$ be a sequence of sequences in $\OO$, $\bsalpha=(\alpha_n)$ a sequence in
$\OO$, $(d\mu^k)$ a sequence in $\frP_0$ and $d\mu\in\frP_0$. If $(\Phi^k_n)$, $(b^k_n)$
are the orthogonal rational functions and inverse iterates associated with $d\mu^k$,
$\bsalpha^k$, and $(\Phi_n)$, $(b_n)$ are the orthogonal rational functions and inverse
iterates associated with $d\mu$, $\bsalpha$, then, for all $n$,
$$
\bsalpha^k \to \bsalpha, \kern5pt d\mu^k \stackrel{*}{\to} d\mu
\kern5pt \Longrightarrow \kern5pt
\bca
b^k_n \stackrel{k}{\rightrightarrows} b_n,
\smallskip \cr
\Phi^k_n \stackrel{k}{\rightrightarrows} \Phi_n, \kern5pt
\Phi^{k*}_n \stackrel{k}{\rightrightarrows} \Phi^*_n \kern5pt
\text{ in } \C\setminus\{\hat\alpha_1,\dots,\hat\alpha_n\}.
\eca
$$

\ep

Now we can prove the strong version of Khrushchev's formula for the orthogonal rational
functions.

\bt[First form of the rational Khruschev's formula] \label{K2}

Let $\bsalpha=(\alpha_n)$ be a sequence in $\OO$, $d\mu\in\frP_0$ and $(\Phi_n)$ the
related orthogonal rational functions. If $(f_n)$ and $(b_n)$ are respectively the
$\bsalpha$-iterates and inverse $\bsalpha$-iterates of $d\mu$, then
$$
f_{\alpha_n}(z;|\Phi_n|^2d\mu) = b_n(z)f_n(z).
$$

\et

\bpr
Let us suppose first that $d\mu(t)=\mu'(t)\,dt$, that is, $d\mu$ is
absolutely continuous. Taking into account that $\cB_{\alpha_n}\cC_{\alpha_n}$ is a
bijection between $\frP_0$ and $\frB$, the fact that $b_nf_n\in\frB$ ensures that
$b_n(z)f_n(z)=f_{\alpha_n}(z;d\sigma_n)$ for some $d\sigma_n\in\frP_0$. In other words,
$$
\frac{1+\zeta_n(z)\,b_n(z)f_n(z)}{1-\zeta_n(z)\,b_n(z)f_n(z)} =
F_{\alpha_n}(z;d\sigma_n).
$$
From (\ref{REF}),
$$
\re\,F_{\alpha_n}(t;d\sigma_n) = \frac{\sigma'_n(t)}{D_R(t,\alpha_n)},
\qquad \text{a.e. } t\in\partial\OO.
$$
On the other hand, Theorem \ref{K1} gives for a.e. $t\in\partial\OO$
$$
\re\left(\frac{1+\zeta_n(t)\,b_n(t)f_n(t)}{1-\zeta_n(t)\,b_n(t)f_n(t)}\right) =
\frac{1-|f_n(t)|^2}{|1-\zeta_n(t)b_n(t)f_n(t)|^2} =
\frac{|\Phi_n(t)|^2}{D_R(t,\alpha_n)}\,\mu'(t).
$$
In consequence, $|\Phi_n|^2\mu'=\sigma'_n$ a.e on $\partial\OO$. Bearing in mind that
$d\sigma_n$ and $d\mu$ are probability measures, the equality
$\int\sigma'_n(t)\,dt=\int|\Phi_n(t)|^2\mu'(t)\,dt=\int|\Phi_n(t)|^2d\mu(t)=1$ shows that
$d\sigma_n$ is absolutely continuous and, thus, $d\sigma_n=|\Phi_n|^2d\mu$. Hence, we
conclude that $b_n(z)f_n(z)=f_{\alpha_n}(z;|\Phi_n|^2d\mu)$.

Consider now an arbitrary measure $d\mu\in\frP_0$, but supported on an infinite subset of
$\partial\OO$. The elements that appear in the rational Khrushchev's formula only depend
on the measure $d\mu$ and the parameters $\alpha_1,\dots,\alpha_n$, but they are independent
of the rest of parameters $\alpha_j$, $j>n$. Therefore, we can suppose without loss of
generality that $B_k\rightrightarrows0$. The absolutely continuous measures
$d\mu^{(k)}=dm_{\alpha_k}/|\Phi_k|^2$ have the same $n$-th orthogonal rational function
as $d\mu$ for $k\geq n$, so
\beq \label{bf-BERNSTEIN}
b_n(z) f^{(k)}_n(z) = f_{\alpha_n}(z;|\Phi_n|^2d\mu^{(k)}), \qquad k \geq n,
\eeq
where $(f^{(k)}_n)_n$ are the $\bsalpha$-iterates of $d\mu^{(k)}$. We know that
$d\mu^{(k)}\stackrel{*}{\to}d\mu$ and $f^{(k)}\rightrightarrows f$ where $f^{(k)}$, $f$
are the S-functions of $d\mu^{(k)}$, $d\mu$ respectively. Hence,
$f^{(k)}_n\stackrel{k}{\rightrightarrows}f_n$ for all $n$ due to Theorem \ref{CONV-gamma}.
Taking the limit $k\to\infty$ in (\ref{bf-BERNSTEIN}), bearing in mind the continuity of
$\cB_{\alpha_n}\cC_{\alpha_n}$, we get Khruschev's formula for $d\mu$.

Finally, suppose that $d\mu\in\frP_0$ is finitely supported. We can obtain $d\mu$ as a
$*$-weak limit of measures $d\mu^k\in\frP_0$ supported on an infinite subset of
$\partial\OO$, for instance, $d\mu^k=\frac{k}{k+1}\,d\mu+\frac{1}{k+1}\,dm$. Denoting
with the superscript ${}^k$ the elements corresponding to the measure $d\mu^k$ and the
sequence $\bsalpha$, we have
\beq \label{bf-INFINITE}
b^k_n(z) f^k_n(z) = f_{\alpha_n}(z;|\Phi^k_n|^2d\mu^k).
\eeq
From the continuity of $\cB_0\cC_0$, Theorem \ref{CONV-gamma} and Proposition \ref{CONV-ORF}
we find that $d\mu^k\stackrel{*}{\underset{k}{\to}} d\mu$ implies that
$\Phi^k_n\stackrel{k}{\rightrightarrows}\Phi_n$ in $\C\setminus\{\hat\alpha_1,\dots,\hat\alpha_n\}$,
$b^k_n\stackrel{k}{\rightrightarrows}b_n$ and $f^k_n\stackrel{k}{\rightrightarrows}f_n$.
Hence, Khrushchev's formula for $d\mu$ is obtained from (\ref{bf-INFINITE}) when $k\to\infty$.
\epr

It could seem surprising that, in the case of a measure with a singular part, the
validity of Khruschev's formula for any sequence $\bsalpha$ is obtained supposing that
the related Blaschke product diverges to 0. Indeed, it is possible to accommodate the
proof of the theorem to a general sequence $\bsalpha$. We simply consider for any fixed
$n$ the new sequence
$\tilde\bsalpha=(\alpha_1,\dots,\alpha_n,\alpha_0,\alpha_0,\alpha_0,\dots)$, so that the
related Blaschke product diverges to 0. Denoting with a tilde the elements related to
$d\mu$ and $\tilde\bsalpha$, we have that $\tilde\Phi_n=\Phi_n$ and $\tilde f_n=f_n$. The
absolutely continuous measures $d\tilde\mu^{(k)}=dm_{\tilde\alpha_k}/|\tilde\Phi_k|^2$
$*$-weak converge to $d\mu$, their $n$-th orthogonal rational functions with respect
$\tilde\bsalpha$ coincide with $\Phi_n$ for $k\geq n$, and their $n$-th
$\tilde\bsalpha$-iterates $\tilde f^k_n$ satisfy $\tilde f^k_n \stackrel{k}{\rightrightarrows}
\tilde f_n=f_n$. So, we get Khrushchev's formula through a limiting process similar
to the one given in the proof of the theorem.

Proposition \ref{f-fa} provides an equivalent version of the strong Khrushchev's formula.

\bc[Second form of the rational Khruschev's formula] \label{K3}

With the notation of Theorem \ref{K2},
$$
f(z;|\Phi_n|^2d\mu) =
-\frac{|\zeta_0(\alpha_n)|}{\zeta_0(\alpha_n)}
\frac{b_n(z)f_n(z)-|\zeta_0(\alpha_n)|}{1-|\zeta_0(\alpha_n)|\,b_n(z)f_n(z)}.
$$

\ec

\section{The indeterminate case} \label{IC}

Proposition \ref{ALL-BERNSTEIN} shows that, in the indeterminate case, to find the limit
points of the sequence of approximants $(R_n/S_n)$ of an S-function $f(z)=f(z;d\mu)$ is
equivalent to find the limit points of the sequence of approximants $(\Psi^*_n/\Phi^*_n)$
of the C-function $F(z)=F(z;d\mu)$ or, alternatively, to find the limit points of the
sequence of measures $(d\mu^{(n)})$. Due to its complexity, the convergence problem of
the $\bsalpha$-continued fraction of $f$ in the indeterminate case will not be completely
addressed in this paper, but we will provide some partial results to understand the
special features of this problem, which does not appear in the polynomial setting. This
discussion will also serve to show an example of application of Khruschev's formula,
whose validity for any sequence $\bsalpha$ makes of it a invaluable tool for studying the
indeterminate case.

Let
$$
\frM_\bsalpha(\bsgamma) = \{d\mu\in\frP_0 : \cS_\bsalpha(d\mu)=\bsgamma\},
\qquad \bsalpha\in\OO, \quad \bsgamma\in\frS.
$$
The determinate case refers to the situation where $\frM_\bsalpha(\bsgamma)$ has only one
measure, otherwise we are in the indeterminate case. The indeterminate case can happen
only if $\bsgamma$ is infinite, so, the measures of $\frM_\bsalpha(\bsgamma)$ are
necessarily infinitely supported in such a situation. Given a sequence $\bsalpha$ in
$\OO$, and bearing in mind the equality $\gamma_n=-z_{n+1}\Lambda_n$, recurrence
(\ref{RR0}) establishes a bijective relation between infinite sequences $\bsgamma$ in
$\D$ and infinite sequences of orthogonal rational functions. Hence, the indeterminate
case corresponds to an infinite sequence of orthogonal rational functions shared by
different measures or, in other words, to an indeterminate rational moment problem:
different measures $d\mu\in\frP_0$ giving the same values of $\int B_n \, d\mu$ for all
$n\in\N$.

The indeterminate rational moment problem was studied in \cite{Al03, BGHN97, BGHN99b,
BGHN99, BGHN07}, following the analysis given in \cite{Ak69} for the polynomial situation
on the real line. In \cite{BGHN97} and \cite[Chapter 10]{BGHN99} it was proved that,
given $\bsalpha=(\alpha_n)$ and $\bsgamma=(\gamma_n)$,
$$
\varDelta(z) = \{F(z;d\mu) : d\mu\in\frM_\bsalpha(\bsgamma)\},
\qquad
z\in\OO_0=\OO\setminus\{\alpha_k\}_{k=0}^\infty,
$$
is always a disk or a point, depending whether we are in the indeterminate case or
not. If $(\Phi_n)$, $(\Psi_n)$ are the orthogonal and second kind rational functions
associated with $\bsalpha$ and $\bsgamma$ throughout a recurrence like (\ref{RR0}), then
$\varDelta(z)$ is a limit of nested disks
\beq
\label{C-DISK}
\varDelta_n(z) =
\{s\in\C : |\Psi^*_n(z)-s\Phi^*_n(z)|\leq|\Psi_n(z)+s\Phi_n(z)|\},
\eeq
which have centers $c_n(z)$ and radius $r_n(z)$ given by
\beq \label{crF}
\kern-7pt
c_n=\frac{\Psi^*_n\overline{\Phi^*_n}+\Psi_n\overline{\Phi_n}}{|\Phi^*_n|^2-|\Phi_n|^2},
\kern20pt
r_n=\frac{|\Psi^*_n\Phi_n+\Psi_n\Phi^*_n|}{|\Phi^*_n|^2-|\Phi_n|^2}
=2{\ts\left|\frac{\varpi_0\,\varpi^*_0}{\varpi_0\kern-1pt(\alpha_0)\,\varpi}\right|}
\frac{|B_{n-1}|}{\sum_{k=0}^{n-1}|\Phi_k|^2},
\eeq
where $\varpi(z)=\varpi_z(z)$ and $(B_n)$ are the Blaschke products related to
$\bsalpha$.

Equivalently, making the substitution $s\to\frac{1+\zeta_0(z)s}{1-\zeta_0(z)s}$ in
(\ref{C-DISK}) and using Proposition \ref{W-ORF}, we find that
$$
\tilde\varDelta(z) = \{f(z;d\mu) : d\mu\in\frM_\bsalpha(\bsgamma)\},
\qquad
z\in\OO_0,
$$
is always a disk or a point, depending whether we are in the indeterminate case or not,
and $\tilde\varDelta(z)$ is a limit of nested disks
\beq \label{S-DISK}
\tilde\varDelta_n(z) =
\{s\in\C : |R_n(z)-sS_n(z)|\leq|S^*_n(z)-sR^*_n(z)|\},
\eeq
with centers $\tilde c_n(z)$ and radius $\tilde r_n(z)$ given by
\beq \label{crf}
\tilde c_n=\frac{R_n\overline{S_n}-S^*_n\overline{R^*_n}}{|S_n|^2-|R^*_n|^2},
\qquad
\tilde r_n=\frac{|S_nS^*_n-R_nR^*_n|}{|S_n|^2-|R^*_n|^2}
=\frac{|B_n|}{\Upsilon_n^{-1}(|S_n|^2-|R^*_n|^2)},
\eeq
where $(R_n)$, $(S_n)$ are the Wall rational functions related to the sequences
$\bsalpha$ and $\bsgamma$ by recurrences (\ref{RR-W2}) and (\ref{RR-W*}).

The determinate case corresponds to the situation where $\varDelta$, or equivalently
$\tilde\varDelta$, is a point in $\OO_0$. In view of the expressions for $r_n$ and
$\tilde r_n$, this occurs iff $B_n$ or $\sum_n|\Phi_n|^2$ diverge in $\OO_0$ (to 0 and
$\infty$ respectively), that is, iff $B_n$ or $\Upsilon_n^{-1}(|S_n|^2-|R^*_n|^2)$
diverge in $\OO_0$ (to 0 and $\infty$ respectively). Therefore, the results of the
previous sections that hold under the divergence of $B_n$, also hold under the divergence
of $\sum_n|\Phi_n|^2$, or equivalently $\Upsilon_n^{-1}(|S_n|^2-|R^*_n|^2)$. For
instance, these conditions ensure the convergence of $(R_n/S_n)$, $(\Psi^*_n/\Phi^*_n)$
and $(d\mu^{(n)})$.

On the contrary, in the indeterminate case, $\varDelta$ and $\tilde\varDelta$ are disks
in $\OO_0$. In this situation $B_n$ necessarily converges, thus $\zeta_n\to1$ and
$\Lim\alpha_n\subseteq\partial\OO$. As for the approximants of the continued fractions,
we only know that $\Lim(\Psi^*_n(z)/\Phi^*_n(z))\subset\varDelta(z)$ and
$\Lim(R_n(z)/S_n(z))\subset\tilde\varDelta(z)$ for any $z\in\OO_0$. However, as we will
see, we can say something more about the limit points of $(\Psi^*_n/\Phi^*_n)$ and
$(R_n/S_n)$ depending on the indeterminate moment problem at hand. Concerning the
possibility of being in the indeterminate case for a given sequence $\bsgamma\in\frS$, we
have the following result.

\bl \label{INDET-ag}

For any infinite sequence $\bsgamma\in\frS$ there exist infinitely many sequences
$\bsalpha$ in $\OO$ such that $\frM_\bsalpha(\bsgamma)$ has more than one measure.

\el

\bpr
Let $\bsgamma\in\frS$ be infinite. We will find sequences $\bsalpha$ in $\OO$ such
that $B_n$ and $\sum_n|\Phi_n|^2$ converge in $\OO$. There, the inequality
$$
|\Phi^*_n|^2 \leq \frac{\varpi_n(\alpha_n)}{\varpi_{n-1}(\alpha_{n-1})}
\frac{1+|\gamma_{n-1}|}{1-|\gamma_{n-1}|}
\left|\frac{\varpi_{n-1}}{\varpi_n}\right|^2 |\Phi^*_{n-1}|^2,
$$
obtained from (\ref{RR0}), proves that
\beq \label{DES-PHI}
|\Phi_n|^2 \leq |\Phi^*_n|^2 \leq
\left|\frac{\varpi_0}{\varpi_n}\right|^2
\frac{\varpi_n(\alpha_n)}{\varpi_0(\alpha_0)}
\prod_{k=0}^{n-1}\frac{1+|\gamma_k|}{1-|\gamma_k|}.
\eeq
Taking into account that
$$
|\varpi_n(z)| \geq \begin{cases} 1-|z|, \cr \im\,z, \end{cases} \qquad z\in\OO,
$$
(\ref{DES-PHI}) shows that the convergence of $\sum_n|\Phi_n|^2$ is a consequence of the
convergence of
$\sum_n\frac{\varpi_n(\alpha_n)}{\varpi_0(\alpha_0)}\prod_{k=0}^{n-1}\frac{1+|\gamma_k|}{1-|\gamma_k|}$.
This last condition also implies the convergence of
$\sum_n\frac{\varpi_n(\alpha_n)}{\varpi_0(\alpha_0)}$, which, bearing in mind (\ref{BP}),
gives the convergence of $B_n$ too. Therefore, it suffices to choose $\bsalpha$ such that
$\sum_n\frac{\varpi_n(\alpha_n)}{\varpi_0(\alpha_0)}\prod_{k=0}^{n-1}\frac{1+|\gamma_k|}{1-|\gamma_k|}$
converges to ensure that $\bsalpha$ and $\bsgamma$ correspond to the indeterminate case.
\epr

The fact that we are in the indeterminate case does not necessarily imply that
$(R_n/S_n)$ is non convergent. For instance, $R_n=0$ and $S_n=1$ if $\gamma_n=0$ for all
$n$. In such a case $R_n/S_n\rightrightarrows0$. However,
$\Upsilon_n^{-1}(|S_n|^2-|R^*_n|^2)=1$ is always convergent, thus we are in the indeterminate
case whenever $B_n$ converges. Nevertheless, this is not the general situation. The
following example shows that $(R_n/S_n)$ can be actually non convergent in the
indeterminate case. Notice that, from (\ref{crf}),
\beq \label{RS-c}
\frac{R_n}{S_n} - \tilde c_n =
- \frac{B_n}{\Upsilon_n^{-1}(|S_n|^2-|R^*_n|^2)} \frac{\overline{R^*_n}}{S_n},
\eeq
thus, in the indeterminate case, $(R_n/S_n)$ converges iff $(\overline{R^*_n}/S_n)$ does
so.

\bex

Let $z\in${\scriptsize$\begin{cases} (0,1) \cr i(1,+\infty) \end{cases}$} \kern-7pt be
fixed. We will choose $\gamma_n\in(-1,1)$ and $\alpha_n\in${\scriptsize$\begin{cases}
(-1,0) \cr i(0,1) \end{cases}$} \kern-7pt so that $z\in\OO_0$ and
$d_n=\overline{R^*_n(z)}/S_n(z)=R^*_n(z)/S_n(z)$ defines a sequence in $(-1,1)$ given by
$$
d_0=\gamma_0; \qquad
d_n=\frac{\zeta_n(z)\,d_{n-1}+\gamma_n}{1+\gamma_n\zeta_n(z)\,d_{n-1}}, \quad n\geq1,
$$
according to (\ref{R*S-schur}). Consider $\varepsilon_n\in(0,1)$ such that
$\sum_n\varepsilon_n<\infty$. Fix $\gamma_0\in(0,1)$ while, for each $n\geq1$, define
$\alpha_n\in${\scriptsize$\begin{cases} (-1,0) \cr i(0,1) \end{cases}$} \kern-7pt by
$$
\frac{\varpi_n(\alpha_n)}{\varpi_0(\alpha_0)} = \varepsilon_n
\prod_{k=0}^{n-1}\frac{1-|\gamma_k|}{1+|\gamma_k|},
$$
and choose $\gamma_n\in(-1,1)$ such that
$$
\begin{cases}
\max\{-\zeta_n(z)\,d_{n-1},\gamma_0\} < \gamma_n < 1 &
\text{ if }  n \text{ is even},
\cr -1 < \gamma_n < \min\{-\zeta_n(z)\,d_{n-1},-\gamma_0\} &
\text{ if } n \text{ is odd}.
\end{cases}
$$
With this choice $d_n>0$ for even $n$ and $d_n<0$ for odd $n$. Besides,
$\sum_n\frac{\varpi_n(\alpha_n)}{\varpi_0(\alpha_0)}
\prod_{k=0}^{n-1}\frac{1+|\gamma_k|}{1-|\gamma_k|}=\sum_n\varepsilon_n$ converges, thus
we are in the indeterminate case, as follows from the proof of Lemma \ref{INDET-ag}. If
$(d_n)$ converges, necessarily $d_n\to0$. In such a case,
$\gamma_n=(1+\gamma_n\zeta_n(z)\,d_{n-1})\,d_n-\zeta_n(z)\,d_{n-1}$ should converge to 0
too, but this is impossible because $|\gamma_n| \geq \gamma_0 > 0$. In consequence,
$(\overline{R^*_n(z)}/S_n(z))$ does not converge, which means that $(R_n(z)/S_n(z))$ is
non convergent because we are in the indeterminate case.
{\huge$\hfill\centerdot$}

\eex

An interesting question is whether the limit points of $(\Psi^*_n/\Phi^*_n)$ are in the
interior $\varDelta^0$ or the frontier $\partial\kern-1pt\varDelta$ of $\varDelta$, which
is equivalent to a similar question concerning $(R_n/S_n)$ and the interior
$\tilde\varDelta^0$ and frontier $\partial\kern-1pt\tilde\varDelta$ of $\tilde\varDelta$.
The reason is that the measures $d\mu\in\frM_\bsalpha(\bsgamma)$ have special features
depending whether $F(z;d\mu)$ lies on $\varDelta^0(z)$ or $\partial\kern-1pt\varDelta(z)$
(a fact which is independent of $z\in\OO_0$, see \cite{BGHN97} and \cite[Chapter
10]{BGHN99}). For example, the condition $F(z;d\mu)\in\partial\varDelta(z)$ for
$z\in\OO_0$, which defines the so called N-extremal measures, characterizes the measures
$d\mu\in\frM_\bsalpha(\bsgamma)$ such that $(\Phi_n)$ is a basis of $L^2(d\mu)$ (see
\cite{BGHN99b} and \cite[Chapter 10]{BGHN99}). Moreover, if the limit points of $\bsalpha$
do not cover $\T$, the map
$$
\mathop{\frM_\bsalpha(\bsgamma) \to \varDelta(z)}
\limits_{\kern20pt \ds d\mu \longrightarrow F(z;d\mu)}
$$
transforms only one measure into each point of $\partial\kern-1pt\varDelta$, while it
transforms infinitely many measures into each point of $\varDelta^0$ (this is a
consequence of the results in \cite{BGHN07}; notice that this property does not appear
correctly written in \cite[Corollary 10.3.2]{BGHN99}).

The modified approximant $(\Psi^*_n-\tau\Psi_n)/(\Phi^*_n+\tau\Phi_n)$ describes
$\partial\varDelta_n$ when $\tau$ runs over $\T$. Hence, given an arbitrary sequence
$(\tau_n)$ in $\T$, the limit points of $(\Psi^*_n-\tau_n\Psi_n)/(\Phi^*_n+\tau_n\Phi_n)$
lie on $\partial\kern-1pt\varDelta$, i.e., they are C-functions of N-extremal measures,
and any C-function of a N-extremal measure can be obtained as a limit of this kind of
modified approximants (see \cite{BGHN97} and \cite[Chapter 10]{BGHN99}). Using the
relation between orthogonal rational functions and Wall rational functions we see that
analogous results hold for the modified approximants
$(R_n-\tau_nS^*_n)/(S_n-\tau_nR^*_n)$ and the S-functions of N-extremal measures. The aim
of the next propositions is to know if something similar happens to the limit points of
$(\Psi^*_n/\Phi^*_n)$ and $(R_n/S_n)$. This is equivalent to analyze the N-extremality of
the limit points of the sequence of measures $(d\mu^{(n)})$.

Our first result concerning the limit points of $(R_n/S_n)$ states that they lie on
$\tilde\varDelta^0$ when $\bsgamma$ converges to zero quickly enough. In what follows, we
will assume that we are in the indeterminate case.

\bp \label{N1}

If $\sum|\gamma_n|<\infty$, the limit points of $(R_n(z)/S_n(z))$ lie on
$\tilde\varDelta^0(z)$ for any $z\in\OO_0$.

\ep

\bpr
Using (\ref{crf}) and (\ref{RS-c}) we find that
$$
\left| \frac{R_n}{S_n} - \tilde c_n \right| =
\left| \frac{R^*_n}{S_n} \right| \tilde r_n.
$$
In consequence, in the indeterminate case, the limit points of $(R_n(z)/S_n(z))$ lie on
$\varDelta^0(z)$ for any $z\in\OO_0$ iff $\Lim(R^*_n/S_n)$ has no constant unimodular
functions. This is also equivalent to $0\notin\Lim(\Upsilon_n/S_n^2)$, as follows from
the identity
$$
\tilde r_n = \frac{\Upsilon_n}{|S_n|^2} \frac{|B_n|}{1-|R^*_n/S_n|^2},
$$
obtained from (\ref{crf}).

From (\ref{RR-W2}) and Proposition \ref{PROP-W}.5 we find that, in $\OO$,
$$
|S_n| \leq \prod_{k=1}^n (1+|\gamma_k|),
$$
thus
$$
\frac{\Upsilon_n}{|S_n|^2} \geq (1-|\gamma_0|^2)
\prod_{k=1}^n\frac{1-|\gamma_k|}{1+|\gamma_k|}.
$$
Hence, $0\notin\Lim(\Upsilon_n/S_n^2)$ if $\prod_n\frac{1-|\gamma_n|}{1+|\gamma_n|}$ does
not diverge to 0, i.e, if $\sum_n|\gamma_n|$ converges.
\epr

The above result does not hold in the general case, as the following proposition shows.

\bp \label{N2}

If $\limsup|\gamma_n|=1$, \kern-2.4pt at least one limit point of $(R_n(z)/S_n(z))$ lies on
$\partial\tilde\varDelta(z)$ for $z\in\OO_0$.

\ep

\bpr
Equivalently, we will prove a similar statement for $(\Psi^*_n(z)/\Phi^*_n(z))$. From
(\ref{crF}) we find that
\beq \label{PP-c}
\left| \frac{\Psi^*_n}{\Phi^*_n} - c_n \right| = |b_n| r_n.
\eeq
Therefore, in the indeterminate case, the limit points of $(\Psi^*_n(z)/\Phi^*_n(z))$ lie
on $\varDelta^0(z)$ for $z\in\OO_0$ iff $\Lim b_n$ has no unimodular constant functions.

Using (\ref{ID}) and Proposition \ref{bn} we get
$$
(1+\gamma_nb_{n+1})(1-\gamma_n\zeta_nb_n)=1-|\gamma_n|^2.
$$
So, if $\limsup|\gamma_n|=1$, then $\liminf(1-|b_{n+1}|)(1-|b_n|)=0$, which gives
$\limsup|b_n|=1$. Hence, $\Lim(\Psi^*_n(z)/\Phi^*_n(z))\nsubseteq\varDelta^0(z)$ for
$z\in\OO_0$. \epr

The next proposition gives a similar result to the previous one, but with a condition for
the sequence $\bsalpha$ instead of $\bsgamma$. Remember that in the indeterminate case
$\Lim\alpha_n\subset\partial\OO$.

\bp \label{N3}

If $\bsalpha$ has a limit point outside of $\supp(d\mu)$, at least one limit point of
$(R_n(z)/S_n(z))$ lies on $\partial\tilde\varDelta(z)$ for $z\in\OO_0$. Furthermore, if
all the limit points of $\bsalpha$ are outside of $\supp(d\mu)$, all the limit points of
$(R_n(z)/S_n(z))$ lie on $\partial\tilde\varDelta(z)$ for $z\in\OO_0$.

\ep

\bpr
Applying Theorem \ref{CONV-borde} to $\bsalpha$ and the sequence of measures
$(|\Phi_n|^2d\mu)$, and taking into account Theorem \ref{K2}, we get in the indeterminate
case
$$
\ba{c}
1 \notin \Lim (b_nf_n)
\kern5pt \Longrightarrow \kern5pt
\Lim(|\Phi_n|^2d\mu) = \{\delta_\tau : \tau\in\Lim\alpha_n\}
\kern5pt \Longrightarrow \kern5pt
\smallskip \cr \kern-28pt
\kern5pt \Longrightarrow \kern5pt
\Lim\alpha_n\subseteq\supp(d\mu).
\ea
$$
Therefore, $\Lim\alpha_n\nsubseteq\supp(d\mu)$ implies $1\in\Lim(b_nf_n)$, so $\Lim b_n$
contains unimodular constant functions. In such a case, following the arguments in the
proof of Theorem \ref{N2} we find that
$\Lim(\Psi^*_n(z)/\Phi^*_n(z))\nsubseteq\varDelta^0(z)$ for $z\in\OO_0$.

Suppose now that $b_nf_n$ does not converge to 1. Then,
$1\notin\Lim_{\!j}\,(b_{n_j}f_{n_j})$ for a subsequence. Theorem \ref{CONV-borde} applied
to $(\alpha_{n_j})_j$ and $(|\Phi_{n_j}|^2d\mu)_j$ gives
$\Lim_{\!j}\,\alpha_{n_j}\subseteq\supp(d\mu)$, so
$\Lim\alpha_n\cap\supp(d\mu)\neq\emptyset$. We conclude that the condition
$\Lim\alpha_n\cap\supp(d\mu)=\emptyset$ implies $b_nf_n\rightrightarrows1$, thus
$|b_n|\to1$, which, bearing in mind (\ref{PP-c}), ensures that
$\Lim(\Psi^*_n(z)/\Phi^*_n(z))\subseteq\partial\varDelta(z)$ for $z\in\OO_0$.
\epr

Theorems \ref{N1}, \ref{N2} and \ref{N3} can be equivalently formulated as statements
about the limit points of $(\Psi^*_n/\Phi^*_n)$ or, alternatively, about the
N-extremality of the limit measures of $(d\mu^{(n)})$. For instance:

\bi

\item[\ref{N1}.] If $\sum|\gamma_n|<\infty$, none of the limit measures of $(d\mu^{(n)})$ is
N-extremal.

\item[\ref{N2}.] If $\limsup|\gamma_n|=1$, at least one limit measure of $(d\mu^{(n)})$ is N-extremal.

\item[\ref{N3}.] If $\bsalpha$ has a limit point outside of $\supp(d\mu)$,
at least one limit measure of $(d\mu^{(n)})$ is N-extremal. Furthermore, if all the limit
points of $\bsalpha$ are outside of $\supp(d\mu)$, all the limit measures of
$(d\mu^{(n)})$ are N-extremal.

\ei

These results are enough to show the variety of possibilities for the convergence of
$(R_n/S_n)$ in the indeterminate case. Besides, Theorem \ref{N3} is obtained as an
application of Khrushchev's formula, showing its interest for the analysis of problems
related to the indeterminate case. Nevertheless, a complete study of the convergence of
$(R_n/S_n)$ in the indeterminate case should address the following problems:

\be

\item A characterization of the cases where $(R_n/S_n)$ is convergent, together with the
description of the corresponding limit.

\item A complete description of the subset of $\tilde\varDelta$ fulfilled by
$\Lim(R_n/S_n)$.

\item A characterization of the limit points of $(R_n/S_n)$.

\ee

\medskip

\noindent{\bf Acknowledgements}

\medskip

This work was partially realized during two stays of the second author at the Norwegian
University of Science and Technology (NTNU) financed respectively by Secretar\'{\i}a de Estado
de Universidades e Investigaci\'{o}n from the Ministry of Education and Science of Spain and
by the Department of Mathematical Sciences of NTNU. The second author wants to express
his gratitude to the Department of Mathematical Sciences of NTNU for the invitations and
the hospitality during both stays. The work of this author was also partly supported by a
research grant from the Ministry of Education and Science of Spain, project code
MTM2005-08648-C02-01, and by Project E-64 of Diputaci\'on General de Arag\'on (Spain).


\begin{thebibliography}{99}

\bibitem{Ak56}
N.I. Akhiezer,
\emph{Theory of Approximation},
Frederic Ungar Publ. Co., New York, 1956.

\bibitem{AkKr62}
N.I. Akhiezer, M.G. Krein,
\emph{Some Questions in the Theory of Moments},
Transl. Math. Monographs, Vol.2, AMS, Providence, RI, 1962;
Russian original, Kharkov, 1938.

\bibitem{Ak69}
N.I. Akhiezer,
\emph{The Classical Moment Problem and Some Related Questions in Analysis},
Oliver and Boyd, Edinburgh, London, 1965; Russian original, Moscow, 1961.

\bibitem{Al03}
A. Almendral, \emph{Nevalinna parametrization of the solutions to some rational moment
problems}, Analysis {\bf 23} (2003), 107--124.

\bibitem{BGHN97}
A. Bultheel, P. Gonz\'{a}lez-Vera, E. Hendriksen, O. Nj{\aa}stad,
\emph{Orthogonal rational functions and nested disks},
J. Approx. Theory {\bf 89} (3) (1997), 344--371.

\bibitem{BGHN99b}
A. Bultheel, P. Gonz\'{a}lez-Vera, E. Hendriksen, O. Nj{\aa}stad,
\emph{A density problem for orthogonal rational functions},
J. Comput. Appl. Math. {\bf 105} (1999), 199--212.

\bibitem{BGHN99}
A. Bultheel, P. Gonz\'{a}lez-Vera, E. Hendriksen, O. Nj{\aa}stad,
\emph{Orthogonal rational functions},
Cambridge Monographs on Applied and Computational Mathematics, 5,
Cambridge University Press, Cambridge, 1999.

\bibitem{BGHN07}
A. Bultheel, P. Gonz\'{a}lez-Vera, E. Hendriksen, O. Nj{\aa}stad,
\emph{An indeterminate moment problem and Carath\'{e}odory functions},
J. Comput. Appl. Math., in press.

\bibitem{Ge44}
Ya.L. Geronimus,
\emph{On polynomials orthogonal on the circle, on trigonometric moment problem,
and on allied Carath\'eodory and Schur functions},
Mat. Sb. {\bf 15} (1944), 99--130. [Russian]

\bibitem{Ge54}
Ya. L. Geronimus,
\emph{Polynomials Orthogonal on a Circle and Their Applications},
Amer. Math. Soc. Translation {\bf 1954} (1954), no. 104, 79 pp.

\bibitem{Ge61}
Ya. L. Geronimus,
\emph{Orthogonal Polynomials: Estimates, Asymptotic Formulas, and Series of
Polynomials Orthogonal on the Unit Circle and on an Interval},
Consultants Bureau, New York, 1961.

\bibitem{JoNjTh89}
W.B. Jones, O. Nj\aa stad, W.J. Thron,
\emph{Moment theory, orthogonal polynomials, quadrature, and continued fractions
associated with the unit circle},
Bull. London Math. Soc. {\bf 21} (1989), 113--152.

\bibitem{Kh01}
S. Khrushchev,
\emph{Schur's algorithm, orthogonal polynomials, and convergence of Wall's continued fractions
in $L\sp 2({\T})$},
J. Approx. Theory {\bf 108} (2001), no. 2, 161--248.

\bibitem{Kh02}
S. Khrushchev,
\emph{Classification theorems for general orthogonal polynomials on the unit circle},
J. Approx. Theory {\bf 116} (2002), no. 2, 268--342.

\bibitem{Si105}
B. Simon,
\emph{Orthogonal Polynomials on the Unit Circle, Part 1: Classical Theory},
AMS Colloquium Series, American Mathematical Society, Providence, RI, 2005.

\bibitem{Si205}
B. Simon,
\emph{Orthogonal Polynomials on the Unit Circle, Part 2: Spectral Theory},
AMS Colloquium Series, American Mathematical Society, Providence, RI, 2005.

\bibitem{Si}
B. Simon,
\emph{CMV matrices: Five years after},
to appear in the Proceedings of the W.D. Evans 65th Birthday Conference,
arXiv:math.SP/0603093, 2006.


\end{thebibliography}
\end{document}